
\input amstex 
\documentstyle{amsppt} 
\magnification=1200 
\voffset -1cm

\input epsf.sty

\NoBlackBoxes
 
\loadbold
 
\TagsOnRight 
\NoRunningHeads 
 
\def\la{\lambda} 
\def\La{\Lambda} 
\def\om{\omega} 
\def\ep{\varepsilon} 
\def\de{\delta} 
 
\def\C{\Bbb C} 
 
\def\Y{\Bbb Y} 
\def\Z{\Bbb Z} 
\def\A{\Bbb A} 
 
\def\Tab{\operatorname{Tab}} 
 
\def\tht{\thetag}

\def\wha{\widehat a} 
\def\whb{\widehat b} 
\def\half{\tfrac12} 
 
\def\sq{{\ssize\square}} 
\def\ssq{{\sssize\square}} 
\def\const{\operatorname{const}} 
\def\T{\Cal T}  
\def\sgn{\operatorname{sgn}} 

\def\bp{\boldkey p}
 
\def\fp#1#2{#1^{\downarrow\,#2}}

\topmatter 
\title Frobenius--Schur functions \endtitle 
\author Grigori Olshanski, Amitai Regev, and Anatoly Vershik \\ 
(with an appendix by Vladimir Ivanov) 
\endauthor 
 
\abstract We introduce and study a new basis in the algebra of 
symmetric functions. The elements of this basis are called the 
Frobenius--Schur 
functions ($FS$-functions, for short).

Our main motivation for studying the $FS$-functions is the fact that 
they enter 
a formula expressing the combinatorial dimension of a skew Young diagram in 
terms of the Frobenius coordinates. This formula plays a key role in the 
asymptotic character theory of the symmetric groups.    
The $FS$-functions are 
inhomogeneous, and their top homogeneous components coincide with the 
conventional Schur functions ($S$-functions, for short). The 
$FS$-functions are best described in the super realization of the 
algebra of symmetric functions. As supersymmetric functions, the 
$FS$-functions can be characterized as a solution to an interpolation 
problem. 

Our main result is a simple determinantal formula for the transition
coefficients between the $FS$- and $S$-functions. We also establish
the $FS$ analogs for a number of basic facts concerning the
$S$-functions: Jacobi--Trudi formula together with its dual form;
combinatorial formula (expression in terms of tableaux); Giambelli
formula and the Sergeev--Pragacz formula.  
 
All these results hold for a large family of bases 
interpolating between the $FS$-functions and the ordinary $S$-functions.  
\endabstract 
 
\toc 
\widestnumber\head{\S7.} 
\head \S0. Introduction \endhead 
\head \S1. Preliminaries on supersymmetric and shifted symmetric 
functions \endhead 
\head \S2. Frobenius--Schur functions \endhead 
\head \S3. Multiparameter Schur functions \endhead 
\head \S4. Combinatorial formula \endhead 
\head \S5. Vanishing property \endhead 
\head \S6. Sergeev--Pragacz formula \endhead 
\head \S7. Transition coefficients \endhead 
\head{} Appendix by V.~Ivanov: Proof of the combinatorial formula for
multiparameter Schur functions \endhead 
\head{} References \endhead 
\endtoc 
 
\endtopmatter 
 
\document 
 
\head \S0. Introduction \endhead 
 
Let $\La$ denote the algebra of symmetric functions. Recall that 
$\La$ is a graded algebra, isomorphic to the algebra of polynomials 
in the power sums $\bp_1,\bp_2,\dots\,$. A natural homogeneous basis in 
$\La$ is formed by the {\it Schur functions\/} (or $S$-functions, for 
short). An $S$-function is denoted as $s_\mu$, where the index $\mu$ 
is a Young diagram.  
 
Let $x=(x_1,x_2,\dots)$ and $y=(y_1,y_2,\dots)$ be two infinite 
collections of indeterminates. We will mainly deal with the {\it 
super realization\/} of $\La$, which is defined by the following 
specialization of the generators $\bp_1,\bp_2,\dots$:  
$$ 
\bp_k\, \to \, \bp_k(x;y)=\sum_{i=1}^\infty x_i^k + 
(-1)^{k-1}\sum_{i=1}^\infty y_i^k\,.   
$$ 
In this realization, each element $F\in\La$ turns into a 
supersymmetric function in $x,y$, which will be denoted as $F(x;y)$. 
In particular, the supersymmetric Schur function will be denoted as 
$s_\mu(x;y)$.  

Any supersymmetric function $F(x;y)$ is separately symmetric in
$x_i$'s and $y_j$'s, and moreover, when a substitution $x_i=-y_j=t$
is made in $F$, the result does not depend on $t$. The latter
property should be viewed as `invariance' under the nonexisting
`supertransposition' $x_i\leftrightarrow y_j$.
 
Let $\mu$ and $\nu$ be Young diagrams. We are interested in the 
combinatorial function  
$$ 
\nu\mapsto\frac{\dim(\mu,\nu)}{\dim\nu}\,,  
$$ 
where $\mu$ is fixed while $\nu$ ranges over the set of Young 
diagrams; $\dim\nu$ is the number of standard tableaux of shape 
$\nu$; $\dim(\mu,\nu)$ is the number of standard tableaux of skew 
shape $\nu/\mu$ provided that $\mu\subseteq\nu$, and 0 
otherwise.   
 
Write $\nu$ in the {\it Frobenius notation\/}:  
$$ 
\nu=(p_1,\dots,p_d\,|\,q_1,\dots,q_d),   
$$ 
where $d=d(\nu)$ (the {\it depth\/} of $\nu$) is the number of 
diagonal squares in $\nu$, $p_i=\nu_i-i$, $q_i=\nu'_i-i$, and $\nu'$ 
stands for the transposed diagram.  Let $m=|\mu|$, $n=|\nu|$, where 
$|\cdot|$ denotes the number of squares of a Young diagram.  Note that the 
function $s_\mu$ is homogeneous of degree $|\mu|$. 
 
Our starting point is the following result:  if $m\le n$ then 
$$ 
\frac{\dim(\mu,\nu)}{\dim\nu} 
=\frac{s_\mu(x(\nu);y(\nu))+\dots} 
{n(n-1)\dots(n-m+1)}\,,  \tag0.1 
$$ 
where $(x(\nu);y(\nu))$ are the so--called {\it modified Frobenius 
coordinates\/} of $\nu$: 
$$ 
(x(\nu);y(\nu))=(p_1+\half,\dots,p_d+\half,0,0,\dots\,;\, 
q_1+\half,\dots,q_d+\half,0,0,\dots),  \tag0.2 
$$ 
and the dots in the numerator denote the remainder term, which is a  
(super)symmetric 
function of degree strictly less than $|\mu|$, evaluated at $(x(\nu);y(\nu))$.

An equivalent form of \tht{0.1} is as follows (\cite{VK}, \cite{W, III.6}):
$$
\frac{\chi^\nu_{\rho\cup 1^{n-m}}}{\dim\nu}=
\frac{\bp_\rho(x(\nu);y(\nu))+\dots}{n(n-1)\dots(n-m+1)}\,, \tag$0.1'$
$$
where $\rho=(\rho_1,\rho_2,\dots)$ is a partition of $m$, $\chi^\nu$ is the 
irreducible character of the symmetric group $S(n)$ which corresponds
to $\nu$, $\chi^\nu_{\rho\cup 1^{n-m}}$ is the character value on the
conjugacy class in 
$S(n)$ indexed by the partition $\rho\cup 1^{n-m}$ of $n$, and 
$\bp_\rho=\bp_{\rho_1}\bp_{\rho_2}\dots\,.$ The equivalence of
\tht{0.1} and \tht{0.1$'$} follows from the relations
$$
\bp_\rho=\sum_\mu\chi^\mu_\rho\,s_\mu\,, \qquad
\chi^\nu_{\rho\cup1^{n-m}}=\sum_\mu\chi^\mu_\rho\dim(\mu,\nu).
$$
 
Formula \tht{$0.1'$} was one of 
the basic ingredients in the asymptotic approach to Thoma's 
description of characters of the infinite symmetric group \cite{T}, see 
\cite{VK}, \cite{KV2}, \cite{W}. Note that 
$\dim\nu$ is a well--known combinatorial function for which  
several nice formulas are known, so that the complexity of \tht{0.1} 
is mainly caused by the {\it relative dimension function\/} 
$\dim(\mu,\nu)$.   

The aim of the present paper is to study in detail the function which
appears in the numerator of the right--hand side of \tht{0.1}. By the
very definition, this is an inhomogeneous supersymmetric function,
and by our convention, we can view it as an element of $\La$. We call it the 
{\it Frobenius--Schur function\/} with index $\mu$ (or the 
$FS$-function, for short) and we denote it as $Fs_\mu$.  
In this notation, \tht{0.1} takes the form 
$$ 
\frac{\dim(\mu,\nu)}{\dim\nu} 
=\frac{Fs_\mu(x(\nu);y(\nu))}{n(n-1)\dots(n-m+1)}\,. \tag0.3 
$$ 

Our main results about the $FS$-functions, and their key properties are as 
follows: 
 
$\bullet$ The $FS$-functions can be characterized in terms of a {\it 
multivariate interpolation problem.} 

Specifically, given a diagram $\mu$, $Fs_\mu$ 
is the only (up to a scalar multiple) supersymmetric function in 
$(x;y)$ that vanishes at $(x;y)=(x(\nu);y(\nu))$ when 
$|\nu|\le|\mu|$, $\nu\ne\mu$, and does not vanish when $\nu=\mu$. 
\footnote{Actually, vanishing holds for all $\nu$ which do not 
contain $\mu$.} 

This characterization follows at once from a correspondence between the $FS$-
functions and the so--called shifted Schur functions (see below), and the 
interpolation properties of the latter functions established in \cite{Ok1}, 
\cite{OO1}. 
 
$\bullet$ As $Fs_\mu$ differs from $s_\mu$ in lower terms only, the 
$FS$-functions form a basis in $\La$. We find explicitly the 
{\it transition coefficients\/} between these two bases $\{s_\mu\}$ and 
$\{Fs_\mu\}$; they are given by simple determinantal expressions. Specifically, 
we have
$$ 
Fs_\mu=\sum_\nu \det[c_{p_i,p'_j}]\,\det[c_{q_i,q'_j}]\, s_\nu, \tag0.4
$$ 
where $p_1,p_2,\dots$, $q_1,q_2,\dots$ denote the Frobenius coordinates of 
$\mu$, and $p'_1,p'_2,\dots$, $q'_1,q'_2,\dots$ denote the Frobenius coordinates 
of $\nu$. The summation in \tht{0.4} is taken over diagrams $\nu$ with the 
length of the diagonal equal to that of $\mu$ (say, $d$). So, the determinants 
are of order $d$. Finally, the coefficients $c_{pp}$ are as follows:
$$ 
c_{pp'}=\cases 
(-1)^{p-p'}e_{p-p'}(\tfrac12,\tfrac32,\dots,\tfrac{2p-1}2), 
& p'\le p,\\ 0, & p'>p, \endcases \tag0.5 
$$ 
for any $p,p'=0,1,2,\dots\,$, where $e_1, e_2,\dots$ stand for the elementary 
symmetric functions. This implies, in particular, that only diagrams $\nu$ 
contained in $\mu$ may contribute. See Theorem 2.6 and the more general Theorem 
7.3.
 
$\bullet$ Consider the {\it duality map\/} $\omega:\La\to\La$ that 
sends $\bp_k$ to $(-1)^{k-1}\bp_k$ (in the super realization, 
$(\omega(F))(x;y)=F(y;x)$ for $F\in\La$). Then (Proposition 2.2)
$$ 
\omega(Fs_\mu)=Fs_{\mu'}\,, 
$$ 
which is similar to the well--known property $\omega(s_\mu)=s_{\mu'}$ of the 
Schur functions. 
 
$\bullet$ For natural $FS$ analogs of the complete homogeneous symmetric 
functions, the elementary symmetric functions, and the hook functions,
$$
Fh_k=Fs_{(k)}, \quad Fe_k=Fs_{(1^k)}, \quad  
Fs_{(k\,|\,l)}=Fs_{(k+1,1^l)}, 
$$
we get counterparts of the well--known generating series. See formulas 
\tht{2.3}, \tht{2.4}, and Theorem 2.3. 

$\bullet$ There exists a wider family of `multiparameter Schur functions' 
which depend on a doubly infinite string $(a_i)_{i\in\Z}$ of complex 
parameters 
and interpolate between the conventional Schur functions 
(case $a_i\equiv0$) and 
the $FS$-functions (case $a_i=i-\half$). See \S3. 

$\bullet$ Recall that the supersymmetric Schur functions in finitely
many variables can be obtained via the {\it Sergeev--Pragacz
formula,\/} which is an analog of the basic formula for the Schur
polynomials, 
$$ 
s_\mu(x_1,\dots,x_n)=\frac{\det[x_i^{\mu_j+n-j}]}{\det[x_i^{n-j}]}\,. 
\tag0.6  
$$ 
In Theorem 6.1 we get a version of the Sergeev--Pragacz formula for
the multiparameter Schur functions (in particular, for the
$FS$-functions).  
 
$\bullet$ Recall the classical combinatorial formula: 
$$ 
s_\mu(x_1,x_2,\dots)=\sum_T\,  
\prod_{(ij)\in\mu}x_{T(i,j)}\,,    
$$ 
summed over the semi--standard tableaux $T$ of the shape $\mu$, and
note that it also has a super analog. We get a version of this
formula for the $FS$-functions (as well as for multiparameter Schur
functions), see Theorem 4.6 and the Appendix.  
 
$\bullet$ The classical Jacobi--Trudi formula and its dual version,  
$$ 
s_\mu=\det[h_{\mu_i-i+j}], \qquad 
s_\mu=\det[e_{\mu'_i-i+j}],  
$$ 
also have counterparts for the $FS$-functions and multiparameter Schur
functions, 
see the formulas just before Theorem 2.3, formula \tht{3.4}, and Corollary 4.9.
 
$\bullet$ Exactly as in the classical Giambelli formula, we have 
$$ 
Fs_\mu=\det[Fs_{(p_i\,|\,q_j)}],   
$$ 
where $p_1,\dots,p_d, q_1,\dots,q_d$ 
are the Frobenius coordinates for $\mu$, $d=d(\mu)$. Moreover, the
same formula 
holds for the multiparameter Schur functions.
 
We now briefly describe the history of formula \tht{0.1} and the 
relationships between our work and \cite{VK}, \cite{W}, \cite{KO},
\cite{OO1}, \cite{Mo}, \cite{L1}, \cite{L2}.  

The formula \tht{0.1} (or rather \tht{$0.1'$}) first appeared in
\cite{VK, \S5} in connection with the asymptotic approach to Thoma's
classification of the characters of the infinite symmetric group.
However, the original proof of \tht{$0.1'$} was not included in
\cite{VK} and remained unpublished. Then a proof was given in 
\cite{W, III.6, Theorem 6} but this work also remained unpublished.   
 
A first published proof of \tht{$0.1'$}) was
given in \cite{KO, Theorem 5}. This paper used the concept of {\it 
shifted symmetric functions\/} (see \S1 for the definition). In 
\cite{KO}, it was shown that there exists a natural isomorphism 
between the algebra $\La^*$ of shifted  
symmetric functions and the algebra of supersymmetric functions 
(hence, the algebra $\La$). Both algebras are viewed as  
algebras of functions on the set of Young diagrams. Under this isomorphism, a 
shifted symmetric function in the row coordinates 
$(\nu_1,\nu_2,\dots)$ of a Young diagram $\nu$ turns into a supersymmetric 
function in the 
modified Frobenius coordinates \tht{0.2}.  
 
The algebra $\La^*$ of shifted symmetric functions was studied in detail in 
\cite{OO1}. In particular, it was proved that 
$$ 
\frac{\dim(\mu,\nu)}{\dim\nu} 
=\frac{s^*_\mu(\nu_1,\nu_2,\dots)}{n(n-1)\dots(n-m+1)}\,, \tag0.7 
$$ 
where $s^*_\mu\in\La^*$ are the so--called {\it shifted Schur 
functions.\/} By \tht{0.3}, this means that under the isomorphism 
$\La^*\to\La$ mentioned above, $Fs_\mu$ is the image of $s^*_\mu$. For 
$s^*_\mu$, a number of explicit expressions is available, see   
\cite{OO1}. This provides us with explicit formulas for 
$\dim(\mu,\nu)$ expressed in terms of the row coordinates. The 
initial purpose of our work was to get similar results in terms of  
the Frobenius coordinates.   
 
Another source of inspiration for our work was Molev's paper
\cite{Mo}. Molev deals with 
supersymmetric polynomials with finitely many indeterminates 
$(x_1,\dots,x_m;$ $y_1,\dots,y_n)$ and studies multiparameter Schur 
polynomials with an arbitrary sequence of parameters $a=(a_i)_{i\in\Z}$. 
He obtains, among other results, generating series for the one--row 
and one--column polynomials, a combinatorial formula, and a 
Sergeev--Pragacz--type formula. However, Molev's approach is not 
quite consistent with our purposes, because his polynomials are not 
stable as $m,n\to\infty$ and, consequently, do not define 
supersymmetric functions in infinitely many indeterminates. Also, his 
formulas are not symmetric with respect to the duality map $\omega$ that 
interchanges the $x$'s and the $y$'s.  
 
We show how to make Molev's construction stable and 
$\omega$-symmetric. The recipe is simple: it suffices to impose 
the restriction $m=n$ and to introduce in the initial definition of 
the polynomials a certain `shift' that depends on $n$. In this way we 
come to a generalization of the $FS$-functions  
depending on $a=(a_i)$. In our approach, we start with Molev's 
generating series, which needs only a slight adaptation. However, our 
versions of the combinatorial formula and the Sergeev--Pragacz--type   
formula differ from Molev's versions in a more substantial way. 

As was pointed out by Lascoux \cite{L1}, \cite{L2}, \cite{OO1,
\S15.6}, the shifted Schur functions can be introduced by means of
the Schubert polynomials.  Possibly, the $FS$ functions can also be
handled using the techniques of the Schubert polynomials.  

As is well known, there exists a deep analogy between the
conventional Schur functions (also called Schur's $S$-functions) and
Schur's $Q$-functions. The latter are related to projective
representations of symmetric groups just in the same way as the former
are related to ordinary representations. Schur's $Q$-functions span
a subalgebra in $\La$, whose elements can be characterized by a
supersymmetry property of another kind, see, e.g., \cite{Ma1, \S
III.8}, \cite{P}. The results of the present paper have counteparts
for Schur's $Q$-functions: these are due to Ivanov \cite{I1},
\cite{I2}. 

For an application of the $FS$ functions, see \cite{BO, \S6}. 

Finally, we would like to notice a remarkable stability of the 
Giambelli formula: in various generalizations of the Schur functions, 
examined in \cite{Ma1}, \cite{OO1}, \cite{Mo}, and the present paper, 
it remains intact while the Jacobi--Trudi identity requires suitable
modification. According to a general theorem due 
to Macdonald \cite{Ma2}, this effect holds under rather wide 
assumptions. The Giambelli formula is substantially 
exploited in our approach.   
 
The paper is organized as follows. In \S1 we collect definitions 
and results from \cite{KO}, \cite{OO1} which are employed in our work. 
In \S2 we introduce the $FS$-functions, then write down the generating 
series and calculate the transition coefficients between the $S$- and 
$FS$-functions. In \S3 we introduce a wider family of `multiparameter 
Schur functions' depending on sequence of parameters $a=(a_i)$ and we
explain their connection with the generalized factorial supersymmetric Schur 
polynomials studied by Molev \cite{Mo}. In \S4 we state the 
combinatorial formula, which expresses the multiparameter Schur 
functions in terms of tableaux. Its proof is given in the 
Appendix, written by Vladimir Ivanov. \S5 is devoted to an analog of  
the Sergeev--Pragacz formula. In \S6 we calculate the transition 
coefficients between multiparameter Schur functions corresponding to 
different sequences of parameters. 

A short exposition of the present paper is given in \cite{ORV}.

{\it Acknowledgment.\/} We would like to thank Alain Lascoux for
discussions and helpful critical remarks. The work was supported by RFBR
grant 98-01--00303 (G.~O.); ISF grant 6629 (A.~R.); CDRF
grant RM1--2244 and RFBR grant 99--01--00098 (A.~V.). G.~O. and A.~V.
are grateful to the Weizmann Institute of Science for hospitality. We
would also like to thank the referee for a number of useful
remarks.
 
\head 1. Preliminaries on supersymmetric and shifted symmetric 
functions  \endhead 
 
In this section, our main references are Macdonald's book \cite{Ma1} 
(symmetric functions in general), \cite{BR} and \cite{Ma1} (supersymmetric 
functions), and \cite{OO1} (shifted symmetric functions). We take $\C$ as the 
base field.
 
The algebra $\La$ of symmetric functions can be initially defined 
as the algebra of polynomials $\C[\bp_1,\bp_2,\dots]$. Then it can 
be realized, in different ways, as an algebra of functions, depending on a 
specialization of the generators $\bp_k$. In the 
conventional realization, the generators $\bp_k$ are 
specialized to the Newton power sums, 
$$ 
\bp_k\mapsto \bp_k(x)=\sum_{i=1}^\infty x_i^k\,, 
$$ 
where $x=(x_1,x_2,\dots)$ is an infinite collection of 
indeterminates. Then elements $f\in\La$ turn into symmetric 
functions in $x$. The term `functions' makes sense if one assumes, 
e.g., that only finitely many of $x_i$'s are different from 0.   
 
Let $\La_n$ be the algebra of symmetric polynomials in 
$x_1,\dots,x_n$. \footnote{To emphasize the difference between the 
case of infinitely many indeterminates and  
that of finitely many indeterminates, we will employ the terms 
`functions' and `polynomials', respectively.} In the conventional 
realization, $\La$ is identified  
with $\varprojlim \La_n$, the projective limit taken in the category of 
graded algebras, where the projection $\La_n\to\La_{n-1}$ is defined as 
the specialization $x_n=0$.  
 
We will mainly deal with another realization of $\La$, which may be 
called the {\it super realization.\/} In this realization, the 
generators $\bp_k$ are viewed as {\it super power sums},  
$$ 
\bp_k\mapsto \bp_k(x;y)=\sum_{i=1}^\infty x_i^k 
+(-1)^{k-1}\sum_{i=1}^\infty y_i^k=\bp_k(x)-\bp_k(-y)\,,  
$$  
where $y=(y_1,y_2,\dots)$ is another collection of indeterminates. 
Note that our definition of the super realization of $\La$ agrees
with that of \cite{BR} but slightly differs from that of \cite{Ma1,
Ex. I.3.23}. 
 
Denote by $\La_{m,n}$ the algebra of polynomials in $x_1,\dots,x_m$, 
$y_1,\dots,y_n$, which are separately symmetric in $x$'s and $y$'s 
and which satisfy the following {\it cancellation property}: for any $i$ and 
$j$, the result of the specialization $x_i=-y_j=t$ does not 
depend on $t$. Such polynomials  
are called {\it supersymmetric,\/} and the (double) projective limit 
algebra $\varprojlim\La_{m,n}$ is called the {\it algebra of 
supersymmetric functions\/} in $x$ and $y$. Here the limit is again 
taken in the category of graded algebras, and the projections 
$\La_{m,n}\to\La_{m-1,n}$ and $\La_{m,n}\to\La_{m,n-1}$ are defined 
by the specializations $x_m=0$ and $y_n=0$, respectively.   
 
By a well--known theorem (see \cite{St}, \cite{P, Theorem 2.11}), the
super realization of $\La$ establishes an isomorphism between the
algebra $\La$ and the algebra $\varprojlim\La_{m,n}$ of
supersymmetric functions, so that we may identify these two algebras.
For this reason, we will not introduce a separate notation for the
algebra of supersymmetric functions.  
 
Given $f\in\La$, we write $f(x;y)$ for the corresponding 
supersymmetric function in $x$ and $y$, and we write  
$f(x_1,\dots,x_m;y_1,\dots,y_n)$ for the supersymmetric polynomial 
that is the image of $f(x;y)$ in $\La_{m,n}$.   
In \S\S4--6 it will be convenient to assume $m=n$ and to think about 
$\La$ as of the projective limit algebra $\varprojlim\La_{n,n}$.   
 
By the cancellation property, the algebra $\La_{n,n}$ of
supersymmetric polynomials is invariant
under any shift of variables of the form 
$$
(x_1,\dots, x_n;y_1,\dots,y_n)\,\mapsto\, 
(x_1+r,\dots,x_n+r;y_1-r,\dots,y_n-r), \qquad r\in\C.
$$
Let $T_r$ denote the corresponding automorphism of $\La_{n,n}$:
$$
(T_rf)(x_1,\dots,x_n;y_1,\dots,y_n)= 
f(x_1+r,\dots,x_n+r;y_1-r,\dots,y_n-r).
$$
Applying again the cancellation property we see that the automorphisms
$T_r:\La_{n,n}\to\La_{n,n}$ and $T_r:\La_{n+1,n+1}\to\La_{n+1,n+1}$ are
compatible with the canonical projection $\La_{n+1,n+1}\to\La_{n,n}$.
Thus, for any $r\in\C$, we get an automorphism $T_r$ of the algebra
$\La=\varprojlim\La_{n,n}$ which informally can be written as follows:
$$ 
(T_rf)(x_1,x_2,\dots;y_1,y_2,\dots)= 
f(x_1+r,x_2+r,\dots;y_1-r,y_2-r,\dots) \qquad 
\forall f\in \La. \tag1.1 
$$ 
Equivalently, $T_r$ can be defined by
$$
T_r: \, \bp_k\,\mapsto\,\sum_{j=1}^k\binom{k}{j}r^{k-j}\bp_j\,.
$$
 
Recall that the elements $h_k$ and $e_k$ (the complete homogeneous symmetric 
functions and the elementary symmetric functions) can be introduced 
through the generating series: 
$$  
1+\sum_{k=1}^\infty h_k t^k  
=\exp\left(\sum_{k=1}^\infty \bp_k\,\frac{t^k}{k}\right) 
=\left(1+\sum_{k=1}^\infty e_k t^k  \right)^{-1}\,, 
$$ 
where $t$ is a formal indeterminate.  
In the super realization,  
$$ 
1+\sum_{k=1}^\infty h_k t^k  \,\mapsto \, 
\prod_{i=1}^\infty \frac{1+y_it}{1-x_it}\,,\qquad 
1+\sum_{k=1}^\infty e_k t^k  \,\mapsto \, 
\prod_{i=1}^\infty \frac{1+x_it}{1-y_it}\,. 
$$ 
 
It will be convenient for us to take $t=\frac1u$ and to redefine the 
generating series for $\{h_k\}$ and $\{e_k\}$ as  
formal series in $\frac 1u$, i.e., as elements of $\La[[\frac1u]]$:  
$$ 
H(u)=1+\sum_{k=1}^\infty \frac{h_k}{u^k}\,, \qquad 
E(u)=1+\sum_{k=1}^\infty \frac{e_k}{u^k}\,.  \tag1.2 
$$ 
Then their super specialization takes the form 
$$ 
H(u)(x;y)=\prod_{i=1}^\infty \frac{1+y_i/u}{1-x_i/u}\,,  \qquad 
E(u)(x;y)=\prod_{i=1}^\infty \frac{1+x_i/u}{1-y_i/u}\,. \tag$1.2'$  
$$ 
Obviously, 
$$ 
H(u)E(-u)=1. \tag1.3 
$$ 
 
\proclaim{Proposition 1.1} In terms of the generating series 
\tht{1.2}, the automorphisms $T_r$ defined in \tht{1.1} 
act on the elements $h_k$ and $e_k$ as follows{\rm:} 
$$ 
T_r(H(u))=H(u-r), \qquad T_r(E(u))=E(u+r).   
$$ 
\endproclaim 
 
It should be noted that for any formal series in $\frac 1u$ (in 
contrast to a series in $u$), the change of a variable  
$u\mapsto u+\operatorname{const}$ makes sense.  
 
\demo{Proof} This follows from the equalities 
$$ 
\frac{1+(y_i-r)/u}{1-(x_i+r)/u}=\frac{1+y_i/(u-r)}{1-x_i/(u-r)}, \qquad 
\frac{1+(x_i+r)/u}{1-(y_i-r)/u}=\frac{1+x_i/(u+r)}{1-y_i/(u+r)}\,. 
$$ 
\qed 
\enddemo 
 
Recall that the {\it duality map\/} is defined as an algebra isomorphism
$\omega: \La\to\La$ such that $\omega(\bp_k)=(-1)^{k-1}\bp_k$. In the super 
realization, $\omega$ reduces to interchanging $x$ and $y$. We have 
$\omega(h_k)=e_k$ and $\omega(e_k)=h_k$.  
 
The {\it Schur function\/} $s_\mu$ indexed by a Young diagram $\mu$ can be  
introduced through the {\it Jacobi--Trudi formula}: 
$$ 
s_\mu=\det[h_{\mu_i-i+j}]\,, 
$$ 
where, by convention, $h_0=1$, $h_{-1}=h_{-2}=\dots=0$, and the order 
of the determinant is any number greater or equal to $\ell(\mu)$, the 
number of nonzero row lengths of $\mu$. The Schur functions form a 
homogeneous basis in $\La$ (it is convenient to agree that 
$s_\varnothing=1$).  
 
Given $m$ and $n$, the supersymmetric Schur polynomial 
$s_\mu(x_1,\dots,x_m;y_1,\dots,y_n)$ in $m+n$ indeterminates does not 
vanish identically if and only if the diagram $\mu$ does not contain 
the square $(m+1,n+1)$. Such polynomials form a basis in $\La_{m,n}$.  
 
We have  
$$ 
\omega(s_\mu)=s_{\mu'}\,, 
$$ 
so that ({\it dual version\/} of Jacobi--Trudi, or {\it 
N\"agelsbach--Kostka formula\/})  
$$ 
s_\mu=\det[e_{\mu'_i-i+j}]\,.  
$$ 
 
For $p,q=0,1,\dots$, let $(p\,|\,q)$ denote the `hook' Young diagram 
$(p+1,1^q)$, and let $s_{(p\,|\,q)}$ denote the corresponding `hook 
Schur function'. \footnote{Note that in \cite{BR}, the term `hook 
Schur functions' has another meaning.} The hook Schur functions form 
a convenient system of generators of $\La$. However, they are 
not algebraically independent. The relations between them are as 
follows:  
$$ 
s_{(p+1\,|\,q)}+s_{(p\,|\,q+1)}=s_{(p\,|\,0)}s_{(0\,|\,q)}\,, 
\qquad p,q=0,1,2,\dots   \tag1.4 
$$ 
(see, e.g., \cite{Ma1, Ex. I.3.9}). This is equivalent to 
$$ 
1+(u+v)\sum_{p,q=0}^\infty \frac{s_{(p\,|\,q)}}{u^{p+1} v^{q+1}} 
=\left(1+\sum_{p=1}^\infty \frac{h_p}{u^p}\right) 
\left(1+\sum_{q=1}^\infty \frac{e_q}{v^q}\right)=H(u)E(v). \tag1.5 
$$ 
 
In the {\it Frobenius notation,\/} a diagram is written as  
$\mu=(p_1,\dots,p_d\,|\,q_1,\dots,q_d\,)$ or, in more detail, 
$$ 
\mu=(p_1(\mu),\dots,p_d(\mu)\,|\,q_1(\mu),\dots,q_d(\mu)), 
$$ 
where $d=d(\mu)$, the {\it depth\/} of $\mu$, is the number of 
diagonal squares, and $p_i=\mu_i-i$, $q_i=\mu'_i-i$.  
 
A very useful expression of the Schur functions is given by the {\it 
Giambelli formula}:  
$$ 
s_\mu=\det[s_{(p_i\,|\,q_j)}]_{i,j=1}^d\,. 
$$ 
 
We proceed to the algebra $\La^*$ of {\it shifted symmetric 
functions.\/} This is a filtered algebra such that the associated 
graded algebra is canonically isomorphic to $\La$. Its definition is 
parallel to that of the algebra $\La$ in its conventional 
realization. First, let $\La^*_n$ be the subalgebra in 
$\C[x_1,\dots,x_n]$ formed by the polynomials which are symmetric in 
`shifted' variables $x'_j=x_j-j$, $j=1,\dots,n$. Define the 
projection map $\La^*_n\to\La^*_{n-1}$ as the specialization $x_n=0$ 
and note that this projection preserves the filtration defined by 
ordinary degree of polynomials. Now  set 
$\La^*=\varprojlim\La^*_n$, where the limit is understood in the 
category of filtered algebras.  
 
The definition of the {\it shifted Schur functions\/} 
$s^*_\mu\in\La^*$ is parallel to the definition of the conventional 
Schur functions via formula \tht{0.6}; the difference is that the 
ordinary powers $x^m$ are replaced by the {\it falling factorial powers} $$ 
\fp xm=\cases x(x-1)\dots(x-m+1), & m\ge1, \\ 
1, & m=0, \endcases 
$$  
and a shift of variables is introduced: 
$$ 
s^*_\mu(x_1,\dots,x_n)= 
\frac{\det[\fp{(x_i+n-i)}{\mu_j+n-j}]}{\det[\fp{(x_i+n-i)}{n-j}]}\,. 
\tag1.6  
$$ 
Clearly, this is a shifted symmetric polynomial of degree  
$|\mu|=\sum_i \mu_i$. Moreover, the formula is stable as $n\to\infty$ 
and, thus, indeed determines an element of $\La^*$. We agree that 
$s^*_\varnothing=1$.  
 
We have  
$$ 
s^*_\mu(x_1,\dots,x_n)=s_\mu(x_1,\dots,x_n)\,+\, 
\text{lower terms}. 
$$ 
It follows that, under the canonical isomorphism 
$\operatorname{gr}\La^*=\La$, the highest term of $s^*_\mu$ coincides 
with $s_\mu$. This implies, in particular, that the shifted Schur 
functions form a basis in $\La^*$. 
 
Recall the definition of the {\it relative dimension function\/} 
$\dim(\mu,\nu)$, where $\mu$ and $\nu$ are arbitrary Young diagrams: 
 
If $\mu\subseteq\nu$ then $\dim(\mu,\nu)$ is the number of standard 
tableaux of skew shape $\nu/\mu$, i.e., the number of chains of 
diagrams of the form 
$$ 
\mu=\la^0\subset\la^1\subset\dots\subset\la^k=\nu,  
\qquad k=|\nu|-|\mu|. 
$$ 
In particular, $\dim(\mu,\mu)=1$. Next, if $\mu$ is not contained in $\nu$ 
then $\dim(\mu,\nu)=0$. Finally, we also set 
$\dim\mu=\dim(\varnothing,\mu)$, which is equal to the number of standard 
tableaux of shape $\mu$. 
 
The relative dimension can be expressed in terms of the shifted Schur 
functions: 
 
\proclaim{Proposition 1.2} For arbitrary Young diagrams 
$\mu,\nu$,  
$$ 
\frac{\dim(\mu,\nu)}{\dim\nu} 
=\frac{s^*_\mu(\nu_1,\nu_2,\dots)}{\fp nm},  
\qquad m=|\mu|, n=|\nu|.  
$$ 
\endproclaim 

Several proofs of this result are given in \cite{OO1, Theorem 8.1}.
The argument presented below is a modification of one of them. 

\demo{Outline of proof} The starting point is the well--known formula
$$
s_\mu\cdot \bp_1=\sum_{\nu:\, \nu\supset\mu, \, |\nu|=|\mu|+1}
s_\nu\,.
$$
It follows that
$$
\bp_1^n=\sum_{\nu:\, |\nu|=n}\dim\nu\cdot s_\nu
$$
and, more generally,
$$
s_\mu\cdot\bp_1^{n-m}=\sum_{\nu:\, |\nu|=n}\dim(\mu,\nu)\, s_\nu\,, 
\qquad |\mu|=m, \quad n\ge m.
$$
Arguing as in \cite{Ma1, Ex. I.1.7} we get from these formulas 
the following ones:
$$
\gather
\dim\nu=n!\,\det\left[\frac1{(\nu_i-i+j)!}\right]_{i,j=1}^l
=\frac{n!\,\prod_{1\le i<j\le l}(\nu_i-\nu_j-i+j)}
{\prod_{i=1}^l(\nu_i+l-i)!}\,,\\
\dim(\mu,\nu)=(n-m)!\,
\det\left[\frac1{(\nu_i-\mu_j-i+j)!}\right]_{i,j=1}^l\,,
\endgather
$$
where $l$ is any sufficiently large natural number ($l\ge n$ is
enough) and $\frac1{k!}=\frac1{\Gamma(k+1)}$ equals 0 when 
$k=-1,-2,\dots$\,. 

This implies
$$
\frac{\dim(\mu,\nu)}{\dim\nu}=
\frac1{\fp nm}\,\cdot\,
\det\left[\frac1{(\nu_i-\mu_j-i+j)!}\right]_{i,j=1}^l\,\cdot\,
\frac{\prod_{i=1}^l(\nu_i+l-i)!}
{\prod_{1\le i<j\le l}(\nu_i-\nu_j-i+j)}\,.
$$
After simple transformations we get
$$
\frac{\dim(\mu,\nu)}{\dim\nu}=
\frac1{\fp nm}\,
\frac{\det\left[\fp{(\nu_i+l-i)}{\mu_j+l-j}\right]_{i,j=1}^l}
{\prod_{1\le i<j\le l}(\nu_i-\nu_j-i+j)}\,,
$$
as required. \qed
\enddemo

\proclaim{Proposition 1.3} For arbitrary Young diagrams $\mu, \nu$, 
$$ 
\gather 
s^*_\mu(\nu_1,\nu_2,\dots)=0 \qquad \text{unless $\mu\subseteq\nu$,} 
\tag1.7\\  
s^*_\mu(\mu_1,\mu_2,\dots)=\frac{|\mu|!}{\dim\mu}\ne0. \tag1.8 
\endgather 
$$ 
\endproclaim 
 
\demo{Proof} See \cite{Ok1} and \cite{OO1, Theorem 3.1}. \qed 
\enddemo 
 
In the same way as for the conventional Schur functions, we set  
$$ 
h^*_k=s^*_{(k)},\qquad e^*_k=s^*_{(1^k)},\qquad 
s^*_{(p\,|\,q)}=s^*_{(p+1,1^q)}.  
$$ 
 
The generating series for $\{h^*_k\}$ and $\{e^*_k\}$ are defined by analogy 
with 
\tht{1.2} but the ordinary powers of $u$ are replaced by the falling 
factorial powers:  
$$ 
H^*(u)=1+\sum_{k=1}^\infty \frac{h^*_k}{\fp uk}\,, \qquad 
E^*(u)=1+\sum_{k=1}^\infty \frac{e^*_k}{\fp uk}\,. 
$$ 
 
\proclaim{Proposition 1.4} The specialization of these series at 
$x=(x_1,x_2,\dots)$ takes the form 
$$ 
H^*(u)(x)=\prod_{i=1}^\infty \frac{1+i/u}{1+(i-x_i)/u}\,,  \qquad 
E^*(u)(x)=\prod_{i=1}^\infty \frac{1+(x_i-i+1)/u}{1+(-i+1)/u}\,, 
$$ 
and the following relation holds {\rm(}cf. \tht{1.3}{\rm)}:  
$$ 
H^*(u)E^*(-u-1)=1. \tag1.9 
$$ 
\endproclaim 
 
\demo{Proof} See \cite{OO1, \S12}. \qed 
\enddemo 
 
Since $\La^*$ is isomorphic to the algebra of polynomials in 
$h^*_1,h^*_2,\dots$, we may define a 1-parameter family $\{T_r\}$, 
$r\in\C$, of automorphisms of the algebra $\La^*$ by  
$$ 
T^*_r(H^*(u))=1+\sum_{k=1}^\infty \frac{T^*_r(h^*_k)}{\fp uk} 
=H^*(u-r).  
$$ 
Together with \tht{1.9} this implies 
$$ 
T^*_r(E^*(u))=1+\sum_{k=1}^\infty \frac{T^*_r(e^*_k)}{\fp uk} 
=E^*(u+r).  
$$ 
It readily follows that 
$$ 
\gather 
T^*_r(h^*_k)=h^*_k\,+\,\{\text{a linear combination of 
$h^*_{k-1},\dots, h^*_1, 1$}\},   \\ 
T^*_r(e^*_k)=e^*_k\,+\,\{\text{a linear combination of 
$e^*_{k-1},\dots, e^*_1, 1$}\}.   
\endgather 
$$ 
 
\proclaim{Proposition 1.5} The following analogs of the Jacobi--Trudi 
formula and its dual version hold{\rm:}  
$$ 
s^*_\mu=\det[T^*_{j-1}(h^*_{\mu_i-i+j})] 
=\det[T^*_{-j+1}(e^*_{\mu'_i-i+j})].  
$$ 
\endproclaim 
 
\demo{Proof} See \cite{OO1, Theorem 13.1}. \qed 
\enddemo 
 
The following general result is due to Macdonald: 
 
\proclaim{Proposition 1.6} Let $A$ 
be a commutative algebra and  
$\{h_{k,r}\}$ be a double sequence of elements of $A$. Here 
$k=1,2,\dots$, $r=0,1,\dots\,.$ We also agree that $h_{0,r}\equiv1$ and 
$h_{k,r}\equiv0$ for all $k<0$. For an arbitrary Young diagram $\mu$, 
set 
$$ 
S_\mu=\det[h_{\mu_i-i+j,j-1}],  
$$ 
where the order of the determinant is any number $\ge\ell(\mu)$. 
Finally, set  
$$ 
S_{(p\,|\,q)}=S_{(p+1,1^q)}. 
$$ 
 
Then the Giambelli formula holds{\rm:} 
$$ 
S_\mu=\det[S_{(p_i\,|\,q_j)}]_{i,j=1}^d\,, 
\qquad \text{where $p_i=p_i(\mu)$, $q_i=q_i(\mu)$, $d=d(\mu)$.}  
$$ 
\endproclaim 
 
\demo{Proof} See \cite{Ma1, Ex. I.3.21} or \cite{Ma2}. \qed 
\enddemo 
 
As a corollary we get the following result, which was pointed out in 
\cite{OO1, Remark 13.2}:  
 
\proclaim{Proposition 1.7} For shifted Schur functions, the Giambelli 
formula remains intact{\rm:} 
$$ 
s^*_\mu=\det[s^*_{(p_i(\mu)\,|\,q_j(\mu))}]. 
$$ 
\endproclaim 
 
\demo{Proof} We employ Proposition 1.5 and apply Proposition 1.6 to 
the two--parameter family $h_{k,r}=T^*_r h_k^*$. \qed 
\enddemo 
 
Following \cite{KO}, we define an algebra isomorphism 
$\varphi:\La^*\to\La$ by  
$$ 
\varphi(H^*(u))=H(u+\half).   
$$ 
Clearly, $\varphi$ intertwines the automorphisms $T_r:\La\to\La$ and 
$T^*_r:\La^*\to\La^*$. The next result yields a  
characterization of $\varphi$. The key idea is to realize both $\La^*$ 
and $\La$ as algebras of functions on Young diagrams; here the 
realization of $\La^*$ is defined in terms of the row coordinates 
while the realization of $\La$ requires Frobenius coordinates. 
 
\proclaim{Proposition 1.8} Let $f\in\La^*$ be arbitrary. For any Young  
diagram $\la$,  
$$ 
f(\la_1,\la_2,\dots)=\varphi(f) 
(p_1+\half,\dots, p_d+\half\,;\, 
\,q_1+\half,\dots,q_d+\half), 
$$ 
where $\la_i$ are the row lengths of $\la$ and  
$(p_1,\dots, p_d\,|\,q_1,\dots,q_d)$ is its Frobenius notation. 
\endproclaim 
 
\demo{Proof} Our argument is a slight simplification of the proof 
given in \cite{KO}. Let $\sq$ range over the squares of $\la$ 
and let $c(\sq)$ denote the content of $\sq$, i.e., $c(\sq)=j-i$ if 
$\sq=(i,j)$. Using Proposition 1.4, we have  
$$ 
\align 
H^*(u)(\la_1,\dots,\la_l)&=\prod_{i=1}^{l}\frac{u+i}{u+i-\la_i}\\ 
&=\prod_{i=1}^{l}\frac{u+i}{u+i-1}\, 
\frac{u+i-1}{u+i-2}\,\dots\, 
\frac{u+i-\la_i+1}{u+i-\la_i}\\ 
&=\prod_{\sq\in\la}\frac{u-c(\sq)+1}{u-c(\sq)}\,.   
\endalign 
$$ 
In the latter product, we may first fix a diagonal hook and then let $\sq$ 
range along this hook. In the $i$th diagonal hook, the content ranges 
from $-q_i$ to $p_i$. From this we conclude that  
$$ 
\gather 
H^*(u)(\la_1,\dots,\la_l)=  
\prod_{i=1}^d\frac{u+q_i+1}{u-p_i} 
=\prod_{i=1}^d\frac{u+\half+(q_i+\half)}{u+\half-(p_i+\half)}\\ 
=H(u+\half)(p_1+\half,\dots,p_d+\half\,; 
q_1+\half,\dots,q_d+\half) 
\endgather 
$$ 
(here we have used \tht{$1.2'$}). \qed 
\enddemo

\head \S2. The Frobenius--Schur functions \endhead 
 
We keep to the notation of \S1.  
 
Define the {\it Frobenius--Schur 
function\/} indexed by a Young diagram $\mu$ as the following element 
of $\La$:  
$$ 
Fs_\mu=\varphi(s^*_\mu). \tag2.1 
$$ 
Here $\varphi$ is the isomorphism $\La^*\to\La$ introduced just before 
Proposition 1.8.

Then the following characterization theorem holds: 
 
\proclaim{Theorem 2.1} $Fs_\mu(x;y)$ is the only 
supersymmetric function such that for any Young diagram 
$\nu=(p_1,\dots, p_d\,|\,q_1,\dots,q_d)$,   
$$ 
\frac{\dim(\mu,\nu)}{\dim\nu} 
=\frac{Fs_\mu(p_1+\half,\dots,p_d+\half,\,; 
\,q_1+\half,\dots,q_d+\half)}{n(n-1)\dots(n-m+1)}\,, \tag2.2 
$$ 
where $m=|\mu|$, $n=|\nu|$. 
\endproclaim 
 
\demo{Proof} Formula \tht{2.2} holds by virtue of Proposition 1.2 and 
Proposition 1.8. The uniqueness claim follows from the fact that a 
supersymmetric function is uniquely defined by its values on the set 
of the modified Frobenius coordinates of Young diagrams. \qed 
\enddemo  

As it was already mentioned in \S0, formula \tht{2.2} is important for the
asymptotic proof of Thoma's theorem \cite{T}. 
Generalizations of Proposition 1.2, Theorem 2.1, and Thoma's theorem
are given in \cite{OO2}, \cite{KOO}, \cite{I1}. 
 
\proclaim{Proposition 2.2} We have  
$$ 
\om(Fs_\mu)=Fs_{\mu'}.  
$$ 
\endproclaim 
 
\demo{Proof} This follows from the fact that the left--hand side 
of \tht{2.2} is invariant under transposition of the diagrams, which
interchanges Frobenius coordinates. \qed 
\enddemo 
 
Set 
$$ 
Fh_k=Fs_{(k)}=\varphi(h^*_k),\qquad 
Fe_k=Fs_{(1^k)}=\varphi(e^*_k). 
$$ 
 
The formulas of \S1 together with the definition \tht{2.1} lead to the 
following formulas for the $FS$ functions. 
 
{\it Generating series\/}: 
$$ 
\gather 
1+\sum_{k=1}^\infty \frac{Fh_k}{(u-\half)\dots(u-\frac{2k-1}2)} 
=H(u)=1+\sum_{k=1}^\infty \frac{h_k}{u^k}\tag2.3\\ 
1+\sum_{k=1}^\infty \frac{Fe_k}{(u-\half)\dots(u-\frac{2k-1}2)} 
=E(u)=1+\sum_{k=1}^\infty \frac{e_k}{u^k}\tag2.4 
\endgather 
$$ 
or, in terms of the super realization of $\La$, 
$$ 
\gather 
1+\sum_{k=1}^\infty \frac{Fh_k(x;y)}{(u-\half)\dots(u-\frac{2k-1}2)} 
=\prod_{i=1}^\infty \frac{1+y_i/u}{1-x_i/u}\\ 
1+\sum_{k=1}^\infty \frac{Fe_k(x;y)}{(u-\half)\dots(u-\frac{2k-1}2)} 
=\prod_{i=1}^\infty \frac{1+x_i/u}{1-y_i/u}\,.
\endgather 
$$ 
 
{\it Jacobi--Trudi and N\"agelsbach--Kostka\/}: 
$$ 
\gather 
Fs_\mu=\det[T_{j-1}(Fh_{\mu_i-i+j})]\\ 
=\det[T_{1-j}(Fe_{\mu'_i-i+j})]. 
\endgather 
$$ 
or, in terms of the super realization, 
$$ 
\gather 
Fs_\mu(x;y)=\det[Fh_{\mu_i-i+j}(x+j-1;y-j+1)]\\ 
=\det[Fe_{\mu'_i-i+j}(x-j+1;y+j-1)],
\endgather 
$$ 
where we abbreviate
$$
x+\const=(x_1+\const,x_2+\const,\dots).
$$
 
{\it Giambelli\/}: 
$$ 
Fs_\mu=\det[Fs_{(p_i\,|\,q_j)}]. 
$$ 
 
The next result is a generating series for the elements $Fs_{(p\,|\,q)}$. 
It will be used in Theorem 2.4 and is of independent interest. 
 
\proclaim{Theorem 2.3} We have 
$$ 
\gather 
1+(u+v)\sum_{p,q=0}^\infty 
\frac{Fs_{(p\,|\,q)}} 
{(u-\tfrac12)\dots(u-\tfrac{2p+1}2)(v-\tfrac12)\dots(v-\tfrac{2q+1}2)}\\ 
=H(u) E(v), 
\tag2.5 
\endgather 
$$ 
cf. \tht{1.5}. Equivalently, 
$$ 
Fs_{(p+1\,|\,q)}+Fs_{(p\,|\,q+1)}+(p+q+1)Fs_{(p\,|\,q)} 
=Fs_{(p\,|\,0)}Fs_{(0\,|\,q)}\,, \qquad p,q=0,1,2,\dots , \tag2.6 
$$ 
cf. \tht{1.4}. 
\endproclaim 
 
\demo{Proof} Since $Fs_{(p\,|\,q)}=\varphi(s^*_{(p\,|\,q)})$ and 
$$ 
H(u)=\varphi(H^*(u-\half)),\qquad  
E(v)=\varphi(E^*(v-\half)), 
$$ 
the relation \tht{2.5} is equivalent to  
$$ 
\multline 
1+(u+v)\sum_{p,q=0}^\infty 
\frac{s^*_{(p\,|\,q)}} 
{(u-\tfrac12)\dots(u-\tfrac{2p+1}2)(v-\tfrac12)\dots(v-\tfrac{2q+1}2)}\\
=H^*(u-\half)E^*(v-\half)\\ 
=\left(1+\sum_{p=0}^\infty 
\frac{s^*_{(p\,|\,0)}} 
{(u-\tfrac12)\dots(u-\tfrac{2p+1}2)}\right) 
\left(1+\sum_{q=0}^\infty 
\frac{s^*_{(0\,|\,q)}} 
{(v-\tfrac12)\dots(v-\tfrac{2q+1}2)}\right). 
\endmultline
$$ 
Writing 
$$ 
u+v=(u-\tfrac{2p+1}2)+(v-\tfrac{2q+1}2)+(p+q+1), 
$$ 
we reduce the above formula to the system of relations 
$$ 
s^*_{(p+1\,|\,q)}+s^*_{(p\,|\,q+1)}+(p+q+1)s^*_{(p\,|\,q)} 
=s^*_{(p\,|\,0)}s^*_{(0\,|\,q)}, \qquad 
p,q=0,1,2,\dots\,. \tag2.7 
$$ 
Note that this system is equivalent to \tht{2.6}. 
 
Thus, it suffices to check \tht{2.7}. This is a particular case of a 
Littlewood--Richardson--type rule for the $s^*$-functions, found by Molev 
and Sagan \cite{MS}. Here is an elementary derivation of \tht{2.7}. 
 
Let us expand the product $s^*_{(p\,|\,0)}s^*_{(0\,|\,q)}$ into 
a linear combination of the shifted Schur functions; such an expansion 
exists, because the shifted Schur functions form a basis in $\La^*$. 
By virtue of \tht{1.4} we know the highest terms of this expansion, so 
we get 
$$ 
s^*_{(p\,|\,0)}s^*_{(0\,|\,q)}= 
s^*_{(p+1\,|\,q)}+s^*_{(p\,|\,q+1)}+ 
\sum_{\nu:\, |\nu|\le p+q+1}c(\nu)s^*_{\nu}\,, 
$$ 
where $c(\nu)$ are certain numerical coefficients. Let $X$ be the set of 
those $\nu$'s that enter this sum with  
nonzero coefficients $c(\nu)$. We claim that $X$ contains only the 
diagram $(p\,|\,q)$; here we will use Proposition 1.3 and the fact that 
$|\nu|\le p+q+1$.  
 
Indeed, let $\la$ be a minimal (with respect to 
inclusion) diagram in $X$ and evaluate both sides at the point  
$\la=(\la_1,\la_2,\dots)$. On the right, the result is nonzero, because 
$s^*_\la(\la)\ne0$ while all other terms on the right have zero 
contributions, because neither a diagram $\nu\ne\la$ from $X$ nor 
$(p+1\,|\,q)$ and $(p\,|\,q+1)$ are contained in $\la$.  So, the result 
of the evaluation on the left is nonzero, too. This implies that $\la$ 
contains both $(p\,|\,0)$ and $(0\,|\,q)$, which is only possible for 
$\la=(p\,|\,q)$.  
 
Thus, our expansion takes the form  
$$ 
s^*_{(p\,|\,0)}s^*_{(0\,|\,q)}= 
s^*_{(p+1\,|\,q)}+s^*_{(p\,|\,q+1)}+ 
\operatorname{const}\,\cdot\, s^*_{(p\,|\,q)}\,, 
$$ 
and it remains to find the constant. Evaluating both sides at 
$\la=(p\,|\,q)$ we get  
$$ 
\const=\frac{s^*_{(p\,|\,0)}(\la)s^*_{(0\,|\,q)}(\la)} 
{s^*_{(p\,|\,q)}(\la)}\,. 
$$  
Applying Proposition 1.2 we get 
$$ 
s^*_{(p\,|\,0)}(\la)= 
\frac{\dim((p\,|\,0),(p\,|\,q))}{\dim(p\,|\,q)}\, 
\fp{(p+q+1)}{p+1}=(p+q+1)p! 
$$ 
Similarly, 
$$ 
s^*_{(0\,|\,q)}(\la)=(p+q+1)q! 
$$ 
Finally, by \tht{1.8},
$$ 
s^*_{(p\,|\,q)}(\la)=\frac{(p+q+1)!}{\dim(p\,|\,q)}= 
(p+q+1)p!q! 
$$ 
This implies $\const=p+q+1$, as was required. \qed  
\enddemo 
 
Our aim is to expand the $FS$ functions in the $S$ functions. 
 
\proclaim{Proposition 2.4} We have  
$$ 
Fs_{(p\,|\,q)}=\sum_{p'=0}^p\sum_{q'=0}^q  
c_{pp'}c_{qq'}s_{(p'\,|\,q')}\,, \tag2.8 
$$ 
where 
$$ 
c_{pp'}=\cases 
(-1)^{p-p'}e_{p-p'}(\tfrac12,\tfrac32,\dots,\tfrac{2p-1}2), & p'\le p,\\ 0, & 
p'>p, \endcases  \tag2.9 
$$ 
for any $p,p'=0,1,2,\dots\,.$
\endproclaim 

(Note that $c_{pp}=1$.)

\demo{Proof} By virtue of \tht{2.5} and \tht{1.5}, 
$$ 
\gather 
(u+v)\sum_{p,q=0}^\infty 
\frac{Fs_{(p\,|\,q)}} 
{(u-\tfrac12)\dots(u-\tfrac{2p+1}2)(v-\tfrac12)\dots(v-\tfrac{2q+1}2)}\\ 
=H(u)E(v)-1=(u+v)\sum_{p,q=0}^\infty 
\frac{s_{(p\,|\,q)}} 
{u^{p+1}v^{q+1}}   
\endgather 
$$ 
Note that $u+v$ can be written as 
$(\frac1u+\frac1v)(\frac1u\cdot\frac1v)^{-1}$ and that the algebra of 
formal power series in the indeterminates $\frac1u$, $\frac1v$ has no zero 
divisors. Consequently, we may divide both sides by 
$u+v$, which gives 
$$ 
\sum_{p,q=0}^\infty 
\frac{Fs_{(p\,|\,q)}} 
{(u-\tfrac12)\dots(u-\tfrac{2p+1}2)(v-\tfrac12)\dots(v-\tfrac{2q+1}2)} 
=\sum_{p,q=0}^\infty 
\frac{s_{(p\,|\,q)}} 
{u^{p+1}v^{q+1}}\,.  \tag2.10 
$$  
 
We need the following Lemma.
 
\proclaim{Lemma 2.5}  
If $a_1,a_2,\dots$ and $b_1,b_2,\dots$ are two number sequences, then 
for any $p'=0,1,\dots$, 
$$ 
\frac1{(u-b_1)\dots(u-b_{p'+1})}=\sum_{p=p'}^\infty  
\frac{h_{p-p'}(b_1,\dots,b_{p'+1}\,;\,-a_1,\dots, -a_p)} 
{(u-a_1)\dots(u-a_{p+1})}\,,  
$$ 
where $h_0=1,h_1, h_2,\dots$ denote the conventional complete homogeneous
functions in the super realization of the algebra $\La$. 
\endproclaim 
 
\demo{Proof of the Lemma} Indeed, we have 
$$ 
\frac1{(u-b_1)\dots(u-b_{p'+1})}=\sum_{p=p'}^\infty  
\frac{d(p,p')} 
{(u-a_1)\dots(u-a_{p+1})}    
$$ 
with certain coefficients $d(p,p')$. It is readily seen that $d(p,p')$ 
is equal to the coefficient of $u^{-1}$ in the expansion of  
$$ 
\frac{(u-a_1)\dots(u-a_p)}{(u-b_1)\dots(u-b_{p'+1})}\,\in\C((\tfrac1u)), $$ 
or, equivalently, to the coefficient of $u^{p'-p}$ in the expansion of  
$$ 
\frac{(u-a_1)\dots(u-a_p)}{(u-b_1)\dots(u-b_{p'+1})u^{p-p'-1}}\, 
\in\C[[\tfrac1u]]. 
$$  
By \tht{$1.2'$}, the latter coefficient is exactly  
$h_{p-p'}(b_1,\dots,b_{p'+1}\,;\,-a_1,\dots, -a_p)$. \qed 
\enddemo 
 
Returning to the proof of the proposition, let us apply Lemma 2.5 to 
$$ 
a_1=\tfrac12,\, a_2=\tfrac32,\, a_3=\tfrac52,\dots, \quad 
b_1=b_2=\dots=0. 
$$ 
Then, using the relation 
$$ 
h_{p-p'}(0,\dots,0\,;\,-\tfrac12,-\tfrac32,\dots,-\tfrac{2p-1}2)=
(-1)^{p-p'}e_{p-p'}(\tfrac12,\tfrac32,\dots,\tfrac{2p-1}2)=c_{pp'},
$$ 
we get 
$$ 
\frac1{u^{p'+1}}=\sum_{p=p'}^\infty 
\frac{c_{pp'}}{(u-\tfrac12)\dots(u-\tfrac{2p+1}2)}\,. 
$$ 
Likewise, we have 
$$ 
\frac1{u^{q'+1}}=\sum_{q=q'}^\infty 
\frac{c_{qq'}}{(u-\tfrac12)\dots(u-\tfrac{2q+1}2)}\,, 
$$ 
whence 
$$ 
\frac1{u^{p'+1}v^{q'+1}}=\sum_{p=p'}^\infty\sum_{q=q'}^\infty 
\frac{c_{pp'}c_{qq'}} 
{(u-\tfrac12)\dots(u-\tfrac{2p+1}2) 
(u-\tfrac12)\dots(u-\tfrac{2q+1}2)}\,.   
$$ 
Substituting this into the right--hand side of \tht{2.10} we get 
$$ 
Fs_{(p\,|\,q)})=\sum_{p'=0}^p\sum_{q'=0}^q     
c_{pp'}\, c_{qq'}\, s_{(p'\,|\,q')}\,,      
$$ 
which concludes the proof of Proposition 2.4. \qed 
\enddemo 
 
\proclaim{Theorem 2.6} {\rm(i)} We have the equality
$$ 
Fs_\mu=\sum_\nu c_{\mu\nu} s_\nu, 
$$ 
summed over diagrams $\nu$ which are contained in $\mu$ and have the same number 
of diagonal squares as $\mu$.  
 
{\rm(ii)} Write $\mu$ and $\nu$ in Frobenius notation 
$$ 
\gathered 
\mu=(p_1,\dots,p_d\,|\,q_1,\dots,q_d),\qquad  
p_1>\dots>p_d\ge0,\quad  
q_1>\dots>q_d\ge0,\\ 
\nu=(p'_1,\dots,p'_d\,|\,q'_1,\dots,q'_d),\qquad 
p'_1>\dots>p'_d\ge0,\quad  
q'_1>\dots>q'_d\ge0. 
\endgathered   
$$ 
Then we have  
$$ 
c_{\mu\nu}=\det[c_{p_i,p'_j}]\,\det[c_{q_i,q'_j}], \tag2.11 
$$ 
where the determinants are of order $d$ and the coefficients
$c_{pp'}$ are defined by \tht{2.9}.  
\endproclaim 
 
\demo{Proof} Both $Fs_\mu$ and $s_\nu$ can be expressed via the 
Giambelli formula,
$$ 
Fs_\mu=\det[Fs_{(p_i\,|\,q_j)}],\qquad 
s_\nu=\det[s_{(p'_i\,|\,q'_j)}]. 
$$ 
Consequently, it suffices to prove the claim of the theorem in the 
simplest case when $\mu$ is a hook diagram, $\mu=(p\,|\,q)$, which 
was done in Proposition 2.4. \qed 
\enddemo

\head \S3. Multiparameter  Schur functions 
\endhead  
 
Let $a=(a_i)_{i\in\Z}$ be an arbitrary sequence 
of complex numbers. All the symmetric functions introduced in this 
section will depend on $a$.  
 
First, we define the multiparameter Schur functions $h_{k;a}$ which are
indexed by one row 
diagrams $(k)$, i.e., certain analogs of the complete homogeneous 
functions. To do this we will employ generating functions which are 
formal series in $u^{-1}$, cf. \tht{2.3}: 
$$ 
1+\sum_{k=1}^\infty \frac{h_{k;a}}{(u-a_1)\dots(u-a_k)} 
=H(u)=1+\sum_{k=1}^\infty\frac{h_k}{u^k}. \tag3.1
$$  
Clearly, we have  
$$ 
h_{k;a}=h_k\,+\,\text{lower terms}. \tag3.2
$$ 
This implies, in particular, that $\{h_{k;a}\}_{k=1,2,\dots}$ is a 
system of algebraically independent generators of $\La$.  
 
We agree that  
$$ 
h_{0;a}=1, \qquad h_{-1;a}=h_{-2;a}=\dots=0. 
$$ 
 
We also need the following notation: for $r\in\Z$, let $\tau^r a$ be the 
result of shifting $a$ by $r$ digits to the left, 
$$ 
(\tau^r a)_i=a_{i+r}\,.\tag3.3
$$ 
 
Now we are in a position to define the {\it multiparameter Schur 
function\/} indexed by an arbitrary Young diagram $\mu$: 
$$ 
s_{\mu;a}=\det[h_{\mu_i-i+j; \tau^{1-j}a}] \tag3.4
$$ 
where the order of the determinant is any number greater or equal to 
$\ell(\mu)$, the number of rows in $\mu$. Clearly, 
$h_{k;a}=s_{(k);a}$.   

Note that our definition \tht{3.4} is a particular case of a very
general concept of {\it multi--Schur functions\/} due to Lascoux
\cite{L1} (see also \cite{L2} and \cite{Ma3}). 
 
\proclaim{Proposition 3.1} Multiparameter Schur functions as defined 
above satisfy the Giambelli formula{\rm:} 
$$ 
s_{\mu;a}=\det[s_{(p_i\,|\,q_j);a}],  
$$ 
where the determinant has order $d=d(\mu)$ and 
$p_1,\dots,p_d; q_1,\dots,q_d$ 
denote the Frobenius coordinates of $\mu$. 
\endproclaim 
 
\demo{Proof} This is immediate from Macdonald's result stated above as 
Proposition 1.6. \qed 
\enddemo 
 
Note that if $a\equiv0$ then $s_{\mu;a}=s_\mu$. For arbitrary $a$, it 
follows from \tht{3.2} and \tht{3.4} that  
$$ 
s_{\mu;a}=s_\mu\,+\,\text{lower terms},  
$$ 
which implies that the elements $s_{\mu;a}$ form a basis in 
$\La$.  
 
\proclaim{Proposition 3.2} If $a_i=i-\half$ then $s_{\mu;a}=Fs_\mu$. 
\endproclaim 
 
\demo{Proof} Both functions can be given by a Jacobi--Trudi--type 
formula, see \S2 and \tht{3.4}. Consequently, it suffices to prove 
that for $a_i=i-\half$ 
$$ 
h_{k;a}=Fh_k, \qquad 
h_{k;\tau^{-r}a}=T_r(Fh_k). 
$$ 
The first equality is immediate from the comparison of \tht{3.1} and 
\tht{2.3}. Let us prove the more general second equality. By \tht{3.1} and 
\tht{3.3},  
$$ 
1+\sum_{k=1}^\infty
\frac{h_{k;\tau^{-r}a}}{(u-a_{1-r})\dots(u-a_{k-r})} =H(u). 
$$ 
By our assumption on $a$, we have $a_{i-r}=a_i-r$. Substitute this 
in the latter expression and then replace $u+r$ by $u$. Then we get  
$$ 
1+\sum_{k=1}^\infty \frac{h_{k;\tau^{-r}a}}{(u-a_1)\dots(u-a_k)} 
=H(u-r). 
$$ 
The left--hand side is equal to
$$ 
1+\sum_{k=1}^\infty \frac{h_{k;\tau^{-r}a}}{(u-\tfrac12)
\dots(u-\tfrac{2k-1}2)}  
$$ 
while the right--hand side, by virtue of Proposition 1.1 and
\tht{2.3}, is equal to 
$$ 
1+\sum_{k=1}^\infty \frac{T_r(Fh_k)}{(u-\tfrac12)\dots(u-\tfrac{2k-1}2)}\,.
$$ 
This proves the second equality. \qed 
\enddemo

Thus, the multiparameter Schur functions interpolate between the 
conventional $S$-functions and the $FS$-functions.

\example{Remark 3.3} Following the general philosophy of symmetric 
functions, one can define multiparameter Schur functions indexed by {\it 
skew\/} diagrams $\la/\mu$ by making use of the canonical 
comultiplication $\Delta:\La\to\La\otimes\La$. We recall that 
$\Delta$ is specified by setting  
$\Delta(\bp_k)=\bp_k\otimes1+1\otimes \bp_k$, or, which is the same, it 
corresponds to splitting the collection  
of the variables into two disjoint parts: $x=x'\sqcup x''$  
(and, in the super case, $y=y'\sqcup y''$). Then $s_{\la/\mu;a}$ is 
defined by 
$$ 
\Delta(s_{\la;a})=\sum_{\mu}s_{\mu;a}\otimes s_{\la/\mu;a}.  
$$ 
 
The result of the Appendix shows that $s_{\la/\mu;a}$ vanishes unless 
$\mu\subseteq\la$, and an analogue of \tht{3.4} holds: 
$$ 
s_{\la/\mu;a}=\det[h_{\la_i-\mu_j-i+j;\tau^{\mu_j-j+1}a}].
$$ 
\endexample 
{}\qed 
 
As was emphasized in \S0, our definition of the multiparameter Schur 
functions was suggested by Molev's work \cite{Mo}. In the rest of this 
section we discuss the connection with \cite{Mo}.  
 
Let $\La_{m,n}$ denote the algebra of supersymmetric polynomials in $m+n$ 
variables \cite{BR}, \footnote{Here we say `supersymmetric 
polynomials' in place of `supersymmetric functions' in order to 
emphasize that one deals with finitely many variables.}  
and let $\Y_{m,n}$ denote the set of Young diagrams not containing 
the square $(m+1,n+1)$. It is well known that the conventional
supersymmetric Schur 
polynomials $s_\mu(x_1,\dots,x_m;y_1,\dots,y_n)$, where $\mu$ ranges 
over $\Y_{m,n}$, form a homogeneous basis in $\La_{m,n}$. The algebra 
of supersymmetric functions can be identified with the projective 
limit of the graded algebras $\La_{m,n}$ as both $m,n$ go to infinity.  
 
In \cite{Mo}, Molev introduced a family of multiparameter supersymmetric 
Schur polynomials, which depend on $a$ and are denoted as 
$$ 
s_\mu(x_1,\dots,x_m/y_1,\dots,y_n|a).  \tag{3.5}
$$ 
These polynomials are inhomogeneous, and their top degree homogeneous 
components coincide with the conventional supersymmetric Schur 
polynomials. This implies that the polynomials \tht{3.5} form a basis in 
$\La_{m,n}$. When 
$a\equiv0$, they reduce to the conventional supersymmetric Schur 
polynomials.   
 
Molev's initial definition is given in terms of a combinatorial 
formula. Next, he writes down generating series for the $h$- and 
$e$-functions and establishes a Jacobi--Trudi--type formula and 
its dual analog. He then gets a version of the Sergeev--Pragacz 
formula. His paper also contains a number of other results which will 
not be discussed here.  
 
Recall that a fundamental property of the conventional supersymmetric 
Schur polynomials is their stability: specializing $x_m=0$ gives the 
supersymmetric Schur polynomial in $(m-1)+n$ variables with the same 
index, and similarly for $y_n=0$. Due to the stability property, one 
can define 
the supersymmetric Schur functions in $\infty+\infty$ variables.  
It is pointed out in \cite{Mo} that the polynomials \tht{3.5} lose the 
stability property. Our observation is that stability is 
recovered if we restrict ourselves to the case $m=n$ and slightly 
modify Molev's definition. The exact correspondence between Molev's 
polynomials and our multiparameter Schur functions is as follows:  
 
\proclaim{Proposition 3.4} For any $\mu$ and any $n$, 
$$ 
s_{\mu;a} (x_1,\dots,x_n;y_1,\dots,y_n)= 
s_\mu(x_1,\dots,x_n/y_1,\dots,y_n|\tau^{-n}a)  
$$ 
\endproclaim 
 
\demo{Proof} Assume first that $\mu=(k)$, where $k=1,2,\dots$, and write 
$h_k(\dots)$ instead of $s_{(k)}(\dots)$. By \cite{Mo, \tht{2.6}} we 
have  
$$ 
1+\sum_{k=1}^\infty  
\frac{h_k(x_1,\dots,x_n/y_1,\dots,y_n|a)} 
{(u-a_{n+1})\dots(u-a_{n+k})} 
=\frac{(u+y_1)\dots(u+y_n)}{(u-x_1)\dots(u-x_n)}\,.
$$ 
Replacing $a$ by $\tau^{-n}a$ we get, by \tht{$1.2'$} and \tht{3.1},
$$ 
\gather
1+\sum_{k=1}^\infty  
\frac{h_k(x_1,\dots,x_n/y_1,\dots,y_n|\tau^{-n}a)} 
{(u-a_{1})\dots(u-a_{k})} 
=\frac{(u+y_1)\dots(u+y_n)}{(u-x_1)\dots(u-x_n)}\\
=H(u)(x_1,\dots,x_n;y_1,\dots,y_n)=1+\sum_{k=1}^\infty  
\frac{h_{k;a}(x_1,\dots,x_n;y_1,\dots,y_n)} 
{(u-a_{1})\dots(u-a_{k})}\,.
\endgather
$$ 
This implies that
$$
h_{k;a}(x_1,\dots,x_n;y_1,\dots,y_n)=h_k(x_1,\dots,x_n/y_1,\dots,y_n|
\tau^{-n}a), \tag3.6
$$
which is our claim for $\mu=(k)$. 
 
For a general $\mu$ we employ Molev's Jacobi--Trudi formula 
\cite{Mo,\tht{3.1}}, which gives 
$$ 
s_\mu(x_1,\dots,x_n/y_1,\dots,y_n|a)= 
\det[h_{\mu_i-i+j}(x_1,\dots,x_n/y_1,\dots,y_n|\tau^{-j+1}a)], 
$$ 
whence 
$$ 
\gather 
s_\mu(x_1,\dots,x_n/y_1,\dots,y_n|\tau^{-n}a)= 
\det[h_{\mu_i-i+j}(x_1,\dots,x_n/y_1,\dots,y_n|\tau^{-n-j+1}a)]\\ 
=\det[h_{\mu_i-i+j;\tau^{1-j}a}(x_1,\dots,x_n;y_1,\dots,y_n)], \quad \text{by 
\tht{3.6}}\\
= s_{\mu;a}(x_1,\dots,x_n;y_1,\dots,y_n), \quad \text{by \tht{3.4}.} 
\endgather 
$$ 
This completes the proof. \qed 
\enddemo 
 
It should be emphasized that this correspondence fails when $m\ne n$.

\head \S4. Combinatorial formula\endhead 
 
We attach to $a$ the `dual' sequence $\wha,$ given by  
$$ 
\wha_i=-a_{-i+1}\,. \tag4.1
$$ 
Note that in the $FS$ case, we have $\wha=a$. 
 
Let $\Z'=\{\dots, -\tfrac32, -\tfrac12, \tfrac12,\tfrac32,\dots\}$ 
stand for the set of proper half--integers. Certain formulas will look more 
symmetric if we agree to label the terms of the sequence $a$ by the 
half--integers. For this reason we introduce the alternative notation  
$$ 
a'=(a'_\ep)_{\ep\in\Z'}, \qquad a'_\ep=a_{\ep+1/2}\,. \tag4.2
$$ 
In this notation, \tht{4.1} takes symmetric form: 
$$ 
(\wha)'_\ep=-a'_{-\ep}. \tag4.3
$$
Note that in the $FS$ case, $a'_\ep=\ep$.

Recall that a {\it horizontal strip\/} is a skew Young diagram 
containing at most one square in each column. Dually, a {\it vertical 
strip\/} contains at most one square in each row. More generally, we 
will deal with skew diagrams $\nu$ of the following kind: there 
exists a skew subdiagram $\nu_1\subseteq\nu$ such that $\nu_1$ is a 
horizontal strip while $\nu/\nu_1$ is a vertical strip (equivalently, 
there exists $\nu_2\subseteq\nu$ which is a vertical strip while 
$\nu/\nu_2$ is a horizontal strip). These are exactly skew diagrams 
$\nu$ containing no $2\times2$ block of squares (equivalently, the 
contents of the squares $\sq\in\nu$ are pairwise distinct). Such a 
diagram is called a {\it skew hook\/} \footnote{Other terms: border 
strip, ribbon, see \cite{Ma1}.} if, in addition, it is connected. 
Thus, a skew diagram with no $2\times2$ block of squares is a disjoint 
union of skew hooks.   
 
To each skew diagram $\nu$ containing no $2\times2$ block of squares 
we attach a polynomial $f_{\nu;a}(u,v)$ in two variables $u,v$, of 
degree $|\nu|$, as follows.   
 
First, assume $\nu$ is a skew hook. Consider the {\it interior\/}
sides of the squares of the shape $\nu$: an interior side is adjacent
to two squares of 
$\nu$; the total number of the interior sides is equal to $|\nu|-1$. 
To each interior side $s$ we attach the coordinates $(\ep,\de)$ of 
its midpoint, \footnote{According to the `English' manner of 
drawing Young diagrams, we assume that the first coordinate axis is 
directed downwards and the second coordinate axis is directed to the 
right.} and we write $s=(\ep,\de)$. Note that one of the  
coordinates is always half--integral while another coordinate is 
integral. Specifically, if $s$ is a vertical side then $\ep\in\Z'$, 
$\de\in\Z$, and the ends of $s$ are the points $(\ep-1/2,\de)$ and 
$(\ep+1/2,\de)$; if $s$ is a horizontal side then $\ep\in\Z$, 
$\de\in\Z'$, and the ends of $s$ are the points $(\ep,\de-1/2)$, 
$(\ep,\de+1/2)$.  
 
For both vertical and horizontal sides, $\de-\ep\in\Z'$. Using the 
notation \tht{4.2}, we set  
$$ 
f_{\nu;a}(u,v)=(u+v) 
\prod\Sb \text{vertical interior}\\ \text{sides $s=(\ep,\de)$ of $\nu$} 
\endSb (u-a'_{\de-\ep}) 
\prod\Sb \text{horizontal interior}\\ \text{sides $s=(\ep,\de)$ of $\nu$} 
\endSb (v+a'_{\de-\ep}). 
\tag4.4
$$ 
 
For instance, if $\nu=(4,2,2)/(1,1)$ (see the figure below) then
there are 6 squares and 5 interior sides with midpoints
$$ 
(\tfrac52,1),\quad(2,\tfrac32),\quad(1,\tfrac32),\quad 
(\tfrac12,2),\quad(\tfrac12,3),
$$ 
and we have  
$$ 
\gather 
f_{\nu;a}(u,v)=(u+v)(u-a'_{-3/2})(u-a'_{3/2})(u-a'_{5/2}) 
(v+a'_{-1/2})(v+a'_{1/2})\\ 
=(u+v)(u-a_{-1})(u-a_2)(u-a_3)(v+a_0)(v+a_1) 
\endgather 
$$ 

(On the figure, the interior sides and
their midpoints are represented by dotted lines and bold dots,
respectively.) 

\smallskip
\midinsert
\vskip 5pt
\centerline{\epsffile{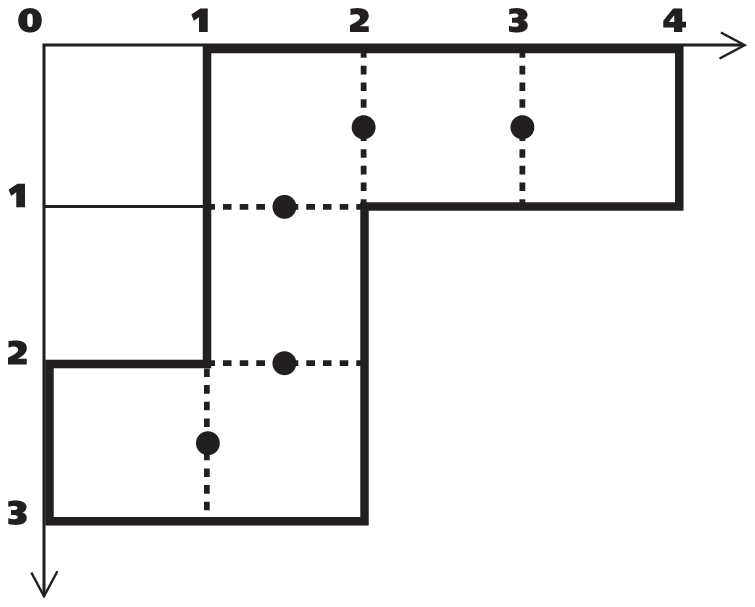}}
\vskip 5pt
\endinsert
\smallskip
 
When $\nu$ is an arbitrary skew diagram with no $2\times2$ block of 
squares, we define $f_{\nu;a}(u,v)$ as the product of the polynomials 
attached to its connected components.  
 
\proclaim{Proposition 4.1} Let $\nu$ be a skew Young diagram 
containing no $2\times2$ block of squares. Then  
$$ 
f_{\nu;a}(v,u)=f_{\nu';\wha}(u,v). 
$$ 
\endproclaim 
 
\demo{Proof} This is immediate from \tht{4.3} and \tht{4.4}. \qed 
\enddemo 
 
Let $\mu$ be a Young diagram. Recall that a 
{\it semistandard\/} (or {\it column--strict\/}) tableau of shape 
$\mu$ is a function $\T(\sq)$ from the squares of $\mu$ to 
$\{1,2,\dots\}$ such that the numbers $\T(\sq)$ weakly increase from 
left to right along the rows and strictly increase down the columns. 
For such a tableau $\T$, the pull--back $\T^{-1}(i)\subset\mu$ is a 
horizontal strip for any $i=1,2,\dots$, see \cite{Ma1, I.5}. 
Dually, for a {\it row--strict\/} tableau, each subset of the form 
$\T^{-1}(i)$ is a vertical strip. Now, we give the following definition: 
 
A {\it diagonal--strict\/} tableau of shape $\mu$ is a function 
$\T(\sq)$ from the squares of $\mu$ to  
$\{1,2,\dots\}$ such that the numbers $\T(\sq)$ weakly increase 
both along the rows (from left to right) and down the columns, and 
strictly increase along the diagonals $j-i=\const$. We will also 
consider diagonal--strict tableaux with entries in $\{1,\dots,n\}$. 
 
Clearly, each subset of form $\T^{-1}(i)$ 
is a skew diagram with no $2\times2$ block of squares, i.e., a 
disjoint union of skew hooks. Thus, a diagonal--strict tableau $\T$ 
with entries in $\{1,\dots,n\}$ may be viewed as a chain of Young 
diagrams,  
$$ 
\varnothing=\mu^{(0)}\subseteq\mu^{(1)}\subseteq\dots 
\subseteq\mu^{(n)}=\mu, 
$$ 
such that $\mu^{(i)}/\mu^{(i-1)}$ has no $2\times2$ block of squares 
for each $i=1,\dots,n$.  
 
Note that the definition of a diagonal--strict tableau also makes 
sense for a skew Young diagram $\mu$.  
 
Given an ordinary or skew Young diagram $\mu$ and 
indeterminates $x=(x_i)$, $y=(y_i)$, consider the combinatorial sum  
$$ 
\Sigma_{\mu;a}(x;y)=\sum_\T 
\prod_{i\ge1}f_{\T^{-1}(i);a}(x_i,y_i)  \tag4.5
$$ 
summed over all diagonal--strict tableaux of shape $\mu$.  
 
By Proposition 4.1, we have 
$$ 
\Sigma_{\mu;a}(y;x)=\Sigma_{\mu';\wha}(x;y).  \tag4.6
$$ 
 
\proclaim{Proposition 4.2} Assume $x_i=y_i=0$ for $i>n$. Then only 
tableaux $\T$ with entries in $\{1,\dots,n\}$ make nonzero 
contributions to the sum \tht{4.5}. 
\endproclaim 
 
\demo{Proof} Indeed, assume $\T$ takes a certain value $i>n$. Then  
for this $i$, the shape $\T^{-1}(i)\subset\mu$ is nonempty. By the 
definition of the polynomials $f_{\nu;a}$ (see \tht{4.4}), 
$f_{\T^{-1}(i);a}(x_i,y_i)$ contains the factor $x_i+y_i$, which is 
zero by the assumption. Consequently, the contribution of $\T$ is 
zero. \qed  
\enddemo 
 
Thus, under the above assumption, the sum \tht{4.5} is actually finite. 
Note that the same holds under the weaker assumption that  
$x_i=-y_i$ for all $i>n$. 
 
\proclaim{Proposition 4.3} The sum \tht{4.5} can also be defined 
by recurrence as follows{\rm:} for any $k<n$, if $\la$ is a skew diagram, 
then  
$$ 
\multline 
\Sigma_{\la;a} 
(x_1,\dots,x_k, x_{k+1},\dots,x_n; y_1,\dots,y_k,y_{k+1},\dots,y_n)\\ 
=\sum_{\mu\subseteq\la} 
\Sigma_{\mu;a}(x_1,\dots,x_k; y_1,\dots,y_k) 
\Sigma_{\la/\mu;a}(x_{k+1},\dots,x_n;y_{k+1},\dots,y_n),  
\endmultline 
$$ 
summed over skew diagrams $\mu$ contained in $\la$, 
and  
$$ 
\Sigma_{\mu;a} (x_1;y_1)=\cases f_{\mu;a}(x_1;y_1), & \text{if $\mu$ 
contains no $2\times2$ block of squares,}\\ 
0, & \text{otherwise.} \endcases 
$$ 
\endproclaim 
 
\demo{Proof} This is evident. \qed 
\enddemo 
 
We proceed with an alternative description of the sum \tht{4.5}. 
Consider the ordered alphabet 
$$ 
\A=\{1'<1<2'<2<\dots\}  
$$ 
and call an {\it $\A$-tableau\/} of shape $\mu$ any map $T(\cdot)$ 
from the set of squares of $\mu$ to the alphabet $\A$ such that: 
 
(*) The symbols $T(\sq)$ weakly increase from left to right 
along each row and down each column. 
 
(**) For each $i=1,2, \dots$, there is at most one symbol $i'$ 
in each row and at most one symbol $i$ in each column. 

This definition (as well as that of diagonal--strict tableaux) is
suggested by the branching rules for the supersymmetric Schur
polynomials, see \cite{BR, \S2}, especially Theorem 2.15 in
\cite{BR}. Note also that the $\A$-tableaux can be obtained via an
appropriate `super' version of the Robinson--Schensted--Knuth
correspondence, see \cite{BR, \S2} and \cite{RS}. Strictly
speaking, the version of the RSK correspondence given in \cite{BR,\S2}
is related to another ordering of the alphabet $\A$. However, the
construction can be readily rephrased to handle our ordering.
Actually, there are many different `super' versions of RSK, related
to different shuffles of the primed and nonprimed indices. This fact
was briefly pointed out at the bottom of page 125 of \cite{BR}. For a
detailed analysis, see \cite{RS}. Note that the RSK correspondence
implies various formulas for the enumeration of the $\A$-tableaux. 
Finally, note that from another point of view, a `super' version of
the RSK correspondence was also discussed in \cite{KV1}.
 
On the other hand, $\A$-tableaux were employed for {\it shifted}
Young diagrams, in the combinatorial formula for the Schur  
$Q$-functions and their factorial analogs, see \cite{Ma1, III.8, 
($8.16'$)}, \cite{I}.   
 
\proclaim{Proposition 4.4} The combinatorial sum \tht{4.5} can be 
written as follows 
$$ 
\Sigma_{\mu;a}(x;y)=\sum_T 
\left(\prod\Sb\ssq\in\mu\\ T(\ssq)=1,2,\dots\endSb 
(x_{T(\ssq)}-a_{c(\ssq)})\,\, 
\prod\Sb\ssq\in\mu\\ T(\ssq)=1',2',\dots\endSb 
(y_{|T(\ssq)|}+a_{c(\ssq)}) 
\right)\,, \tag4.7
$$ 
summed over all $\A$-tableaux of shape $\mu$, where we use the notation 
$|i'|=i$ for $i=1,2,\dots$ and $c(\sq)$  
denotes the content of $\sq$, i.e., if  $\sq=(i,j)$ then $c(\sq)=j-i$. 
\endproclaim 
 
\demo{Proof} Let $\nu$ be a skew hook and $i\in\{1,2,\dots\}$ be 
fixed. Clearly, there exists exactly one diagonal--strict tableau 
$\T$ of shape $\nu$, with entries in $\{i\}$, and we claim that 
there exist exactly two $\A$-tableaux $T$ of the same shape, with 
entries in $\{i,i'\}$. Indeed, let $\sq_1,\dots,\sq_k$ be the 
squares of $\nu$ written down in the order of increasing contents. Then 
$T(\sq_1)$ may be chosen arbitrarily, while for any $r=2,\dots,k$, 
the value of $T(\sq_r)$ is defined uniquely, according to whether 
the squares $\sq_{r-1}$, $\sq_r$ lie in the same row or in the same 
column: in the former case, $T(\sq_r)=i$, and in the latter case, 
$T(\sq_r)=i'$.  
 
Next, if $\nu$ is a skew diagram with no $2\times2$ block of squares 
then the same is true, with the only exception that the number of the 
$T$'s is equal to $2^l$, where $l$ stands for the number of connected 
components of $\nu$.  
 
Now, let $\mu$ be an arbitrary (skew) diagram. To any $\A$-tableau 
$T$ of shape $\mu$ we assign a diagonal--strict tableau $\T$ by 
replacing each primed index $i'$ by the corresponding nonprimed index 
$i$. Conversely, any diagonal--strict tableau $\T$ of shape $\mu$ can 
be obtained in this way from a certain (nonunique) $\A$-tableau $T$. 
To get all such $T$'s, we have to choose, for any nonempty diagram  
$\nu=\T^{-1}(i)$, an arbitrary $\A$-tableau of shape $\nu$, with entries in 
$\{i,i'\}$, as described above.  
 
This means that the right--hand side in \tht{4.7} can be written as a double 
sum, 
where the exterior sum is taken over the diagonal--strict tableaux 
$\T$ of shape $\mu$, and each interior sum is taken over all 
$\A$-tableaux $T$ `over' a fixed $\T$. It follows that the claim of 
the proposition can be reduced to the following one: let $\nu$ by a skew
Young diagram with no $2\times2$ block of squares; then  
$$ 
f_{\nu;a}(u;v)=\sum_{T:\nu\to\{1',1\}} 
\left(\prod\Sb\ssq\in\nu\\ T(\ssq)=1\endSb 
(u-a_{c(\ssq)})\,\, 
\prod\Sb\ssq\in\nu\\ T(\ssq)=1'\endSb 
(v+a_{c(\ssq)}) 
\right)\,, 
$$ 
summed over all $\A$-tableaux of shape $\nu$, with entries in 
$\{1,1'\}$.   
 
Finally, without loss of generality, we may assume that $\nu$ is a 
skew hook. Then, as was shown above, there are exactly two $T$'s, so 
that the sum in the right--hand side of the last formula consists of two 
summands. On the 
other hand, the left--hand side is given by \tht{4.4}. Writing in that 
expression the factor  
$u+v$ as the sum of $u-a_c$ and $v+a_c$, where $c=c(\ssq_1)$ is the 
smallest content, we split the left--hand side into two summands, too. Then the 
desired equality is readily verified. \qed  
\enddemo 
 
The following claim is a slight refinement of Proposition 4.4: 
 
\proclaim{Proposition 4.5} Assume $x_i=y_i=0$ {\rm(}or, more generally, 
$x_i=-y_i${\rm)} for all $i>n$. Then, in the right--hand side of \tht{4.7}, one 
can 
take only tableaux $T$ with entries in $\A_n=\{1'<1<\dots<n'<n\}$. 
\endproclaim 
 
\demo{Proof} Indeed, this follows from Proposition 4.3 and the proof 
of Proposition 4.4. \qed 
\enddemo 
 
Note that the $\A$-tableaux $T$, in contrast to the diagonal--strict 
tableaux $\T$, are not consistent with transposition. In particular, 
the symmetry \tht{4.6} is not evident from \tht{4.7}. However, in certain 
circumstances, it is more convenient to use formula \tht{4.7} than 
formulas \tht{4.4} and \tht{4.5}. 
 
\proclaim{Theorem 4.6 {\rm(Combinatorial formula)}} We have 
$$ 
s_{\mu;a}(x;y)=\Sigma_{\mu;a}(x;y) 
$$ 
where the right--hand side is given by \tht{4.5} or \tht{4.7}. 
\endproclaim 
 
\demo{Proof} See the Appendix. \qed 
\enddemo 

Note that for the first time, an `inhomogeneous' combinatorial
formula probably appeared in \cite{BL1}, \cite{BL2}, see also
\cite{CL}. Other examples can be found in \cite{GG},
\cite{Ma2}, \cite{Mo}, \cite{Ok1}, \cite{OO1}, \cite{Ok2}. See also
further references in \cite{Ok2} to works by Knop, Okounkov, and Sahi
about combinatorial formulas for interpolation Jack and Macdonald
polynomials.  

\proclaim{Corollary 4.7 {\rm(Duality)}} We have 
$$ 
\om(s_{\mu;a})=s_{\mu';\wha}.   
$$ 
\endproclaim 
 
\demo{Proof} Indeed, this follows from Theorem 4.6 and \tht{4.6}. \qed 
\enddemo 
 
Set  
$$ 
e_{k;a}=s_{(1^k);a}\,. 
$$ 
By Corollary 4.7, 
$$ 
\om(h_{k;a})=e_{k;\wha}. \tag4.8
$$ 
 
\proclaim{Corollary 4.8 {\rm(Generating series for $e$-functions)}} We 
have  
$$ 
1+\sum_{k=1}^\infty \frac{e_{k;a}}{(u-\wha_1)\dots(u-\wha_k)} 
=E(u)=1+\sum_{k=1}^\infty\frac{e_k}{u^k}. \tag4.9
$$ 
\endproclaim 
 
\demo{Proof} Applying $\om$ to both sides of \tht{3.1} 
we get
$$
1+\sum_{k=1}^\infty \frac{e_{k;\wha}}{(u-a_1)\dots(u-a_k)} 
=\om(H(u))=E(u). 
$$ 
Next, replacing $a$ by $\wha$, we get \tht{4.9}. \qed 
\enddemo 
 
\proclaim{Corollary 4.9 {\rm(N\"agelsbach--Kostka formula)}} We have 
$$ 
s_{\mu;a}=\det[e_{\mu'_i-i+j;\tau^{j-1}a}] 
$$ 
with the understanding that  
$$ 
e_{0;a}=1, \qquad e_{-1;a}=e_{-2;a}=\dots=0. 
$$ 
\endproclaim 
 
\demo{Proof} This follows from \tht{3.4}, Corollary 4.7, \tht{4.8}, and the fact 
that 
$(\tau^r a)\widehat{\phantom{a}}=\tau^{-r}\wha$. \qed 
\enddemo

\example{Remark 4.10} Let us specialize $a_i=i-\tfrac12$. That is, take 
$a'_{\de-\ep}=\de-\ep$ in \tht{4.4}, and $a_{c(\sq)}=c(\sq)-\tfrac12$ in 
\tht{4.7}. Then Theorem 4.6 turns into a combinatorial formula for the $FS$-
functions. 
\endexample

\head \S5. Vanishing property\endhead 
 
We fix a Young diagram $\mu$ and write it in the Frobenius notation,  
$$ 
\mu=(p_1,\dots,p_d\,|\,q_1,\dots,q_d). 
$$ 
Let $\la$ be an arbitrary diagram, 
$$ 
\la=(P_1,\dots,P_D\,|\,Q_1,\dots,Q_D). 
$$ 
Then $\mu\subseteq\la$ means that $d\le D$, $p_i\le P_i$, $q_i\le 
Q_i$ for $i=1,\dots, d$.  
 
We define a collection of variables $(x(\la);y(\la))$ as follows 
$$ 
\gather 
x(\la)_i=a_{P_i+1},\quad  
y(\la)_i=\wha_{Q_i+1},\quad 1\le i\le D,\\  
x(\la)_i=y(\la)_i=0, \quad i>D. 
\endgather 
$$ 
 
For instance, if $a_i=i-\half$ then $(x(\la),y(\la))$ is exactly the 
collection of the modified Frobenius coordinates of $\la$.  
 
\proclaim{Theorem 5.1 {\rm(Vanishing Theorem)}} If 
$\mu\not\subseteq\la$ then $s_{\mu;a}(x(\la);y(\la))=0$. 
\endproclaim 
 
\demo{Proof} We employ, in a slightly modified 
form, an argument due to Okounkov, cf. \cite{Ok1, proof of Prop. 3.8} 
and \cite{OO1, second proof of Theorem 11.1}.  
 
Assume that $s_{\mu;a}(x(\la);y(\la))\ne0$ and let us prove that 
$\la\supseteq\mu$.  
 
{\it Step 1.\/} Let us prove that $D\ge d$. Employ for 
$s_{\mu;a}(x(\la);y(\la))$ the expression given by Theorem 4.6 and
formulas \tht{4.4} and \tht{4.5}. By Proposition 4.2, we can take in
\tht{4.5} only tableaux $\T$ with entries in $\{1,\dots,D\}$. On the
other hand, the main diagonal in $\mu$ has length $d$ and, by the
definition of a diagonal--strict tableau, it is filled by strictly
increasing numbers. Consequently, $D\ge d$.  
 
{\it Step 2.\/} By Corollary 4.7, the quantity $s_{\mu;a}(x(\la);y(\la))$ does 
not change under $\mu\mapsto\mu'$, $\la\mapsto\la'$, $a\mapsto\wha$. 
Consequently, to conclude that 
$\la\supseteq\mu$, it suffices to prove that  
$$ 
P_1\ge p_1, \dots, P_d\ge p_d\,.  \tag5.1
$$ 
Write the coordinates $x_i(\la)$, $y_i(\la)$, $i=1,\dots,D$, in 
the reverse order, 
$$ 
\gather 
\bar x=(\bar x_1,\dots,\bar x_D)= 
(a_{P_D+1},\dots,a_{P_1+1}), \\ 
\bar y=(\bar y_1,\dots,\bar y_D)= 
(\wha_{Q_D+1},\dots,\wha_{Q_1+1}). 
\endgather 
$$ 
Since $s_{\mu;a}(x;y)$ is symmetric in $x$ and in $y$, we get 
$$ 
s_{\mu;a}(\bar x;\bar y)=s_{\mu;a}(x(\la);y(\la))\ne0. 
$$ 
By Theorem 4.6, $s_{\mu;a}(\bar x;\bar y)$ is given by formula  \tht{4.7}.  Let 
us fix a tableau $T$ which has nonzero contribution to the sum \tht{4.7}. For 
this $T$, we get, in particular, 
$$ 
\prod\Sb\ssq\in\mu\\ T(\ssq)=1,2,\dots\endSb 
(\bar x_{T(\ssq)}-a_{c(\ssq)})\ne0.  \tag5.2
$$ 
We aim to prove that \tht{5.2} implies \tht{5.1}. Note that, by 
Proposition 4.5, $|T(\ssq)|$ takes values in $\{1,\dots,D\}$. 
 
Introduce the notation 
$$ 
(k(1),\dots,k(D))=(P_D+1,\dots,P_1+1),  
$$  
so that 
$$ 
1\le k(1)<\dots<k(D).  \tag5.3
$$ 
In this notation, 
$$ 
\bar x_r=a_{k(r)}, \qquad r=1,\dots,D.  
$$ 
 
Our argument will employ the following evident fact: 

(*) If, for a certain square $\sq\in\mu$, $T(\sq)$ is nonprimed 
then $k(T(\sq))\ne c(\sq)$.

(Indeed, otherwise we would get  
$$ 
\bar x_{T(\ssq)}-a_{c(\ssq)}=a_{k(T(\ssq))}-a_{c(\ssq)}=0, 
$$ 
in contradiction with \tht{5.2}.) 
 
{\it Step 3.\/} Let $\T$ be the diagonal--strict tableau 
corresponding to $T$. By its definition, $\T(\sq)=|T(\sq)|$ for any 
square $\sq\in\mu$. Consider the squares $\sq=(1,j)$ of the first 
row in $\mu$. For these squares, 
$$ 
\T(1,1)\le\T(1,2)\le\dots\le\T(1,\mu_1). 
$$ 
We claim that  
$$ 
k(\T(1,j))\ge j,  \qquad j=1,\dots,\mu_1.  \tag5.4
$$ 
 
Indeed, \tht{5.4} is trivial for $j=1$. Assuming that \tht{5.4} is true for  
$j\le j_0$, let us check it for $j=j_0+1$.  
 
The numbers $k(\T(1,j))$ weakly increase. Therefore, if 
$k(\T(1,j_0+1))\le j_0$ then, by the assumption, 
$$ 
k(\T(1,j_0))=k(\T(1,j_0+1))=j_0. 
$$ 
It follows that $T(1,j_0+1)$ is nonprimed (indeed, if $T(1,j_0+1)$ 
were primed then $T(1,j_0)$ would be the same primed index, which contradicts 
to the 
definition of $\A$-tableaux). Then for the square $\sq=(1,j_0+1)$ we 
get: $T(\sq)$ is nonprimed and $k(T(\sq))=j_0=c(\sq)$, which is in 
contradiction with (*). Thus, we have proved \tht{5.4} by induction.  
 
{\it Step 4.\/} Recall that the numbers $\T(\sq)=|T(\sq)|$ take 
values in $\{1,\dots,D\}$. On the other hand, these numbers strictly 
increase as $\sq$ ranges over any diagonal in $\mu$ from top to
bottom. Together with 
the inequalities \tht{5.3} this implies that the numbers $k(\T(\sq))$ 
also strictly increase along diagonals. Consequently, for any square 
$\ssq=(i,j)\in\mu$ with $i\le j$, 
$$ 
k(\T(i,j))>k(\T(i-1,j-1))>\dots>k(\T(1,j-i+1))\ge j-i+1, \tag5.5
$$ 
where, on the last step, we have used \tht{5.4}.  
 
Let us fix $i=1,\dots,d$ and set $j=\mu_i$, which 
corresponds to the last square in the $i$th row of $\mu$ (since $i\le 
d$, the assumption $j\ge i$ is satisfied). Then $j-i+1=p_i+1$. On the other 
hand, by the definition, each number $k(\,\cdot\,)$ is equal 
to a certain $P_r+1$, and recall that $P_1>\dots>P_D$, so that 
$P_1+1>\dots>P_D+1$. Thus, \tht{5.5} means that, for certain
$r_1<\dots<r_i$, 
$$
P_{r_1}+1>\dots>P_{r_i}+1\ge p_i+1.
$$
Since $r_i\ge i$, it follows that $P_i\ge p_i$, which concludes the 
proof. \qed 
\enddemo 
 
\proclaim{Theorem 5.2} Let $\mu$ be a Young diagram and let the symbols 
$x(\cdot)$, $y(\cdot)$ be as defined in the beginning of the section. We have 
$$ 
s_{\mu;a}(x(\mu);y(\mu))= 
\prod_{(i,j)\in\mu}(a_{\mu_i-i+1}-a_{j-\mu'_j}). \tag5.6
$$ 
\endproclaim 
 
Note that  
$$ 
a_{\mu_i-i+1}-a_{j-\mu'_j}=a'_{\mu_i-i+1/2}-a'_{-\mu'_j+j-1/2}= 
a_{\mu_i-i+1}+\wha_{\mu'_j-j+1}\,, 
$$ 
where $a'=(a'_\ep)$ was defined in \tht{4.2}. Hence, \tht{5.6} can 
be rewritten as follows: 
$$ 
\multline s_{\mu;a}(x(\mu);y(\mu))\\
=\prod_{(i,j)\in\mu}(a'_{\mu_i-i+1/2}-a'_{-\mu'_j+j-1/2})= 
\prod_{(i,j)\in\mu}(a_{\mu_i-i+1}+\wha_{\mu'_j-j+1}). 
\endmultline \tag5.7
$$ 
These two expressions are symmetric with respect to 
$\mu\leftrightarrow\mu'$. One more expression is given below in 
Proposition 5.4.
 
\demo{Proof} As before, write $\mu$ in the Frobenius notation,  
$\mu=(p_1,\dots,p_d\,|\,q_1,\dots,q_d)$. Let $\mu^*\subset\mu$ denote 
the shape which is obtained from $\mu$ by removing the border skew 
hook --- the set of the squares $(i,j)\in\mu$ such that 
$(i+1,j+1)\notin\mu$. In the Frobenius notation, $\mu^*$ is obtained by 
removing the Frobenius coordinates $p_1$ and $q_1$, 
$$ 
\mu^*=(p_2,\dots,p_d\,|\,q_2,\dots,q_d). 
$$ 
Note that $\mu^*$ is the smallest subdiagram in $\mu$ such that 
$\mu/\mu^*$ has no $2\times2$ block of squares.  
 
Let $\bar x(\mu)$ and $\bar y(\mu)$ be obtained from $x(\mu)$ and 
$y(\mu)$ by writing the coordinates in the reversed order, 
$$ 
\gather 
\bar x(\mu)=(\bar x_1,\dots,\bar x_d)= 
(a_{p_d+1},\dots,a_{p_1+1}), \\ 
\bar y(\mu)=(\bar y_1,\dots,\bar y_d)= 
(\wha_{q_d+1},\dots,\wha_{q_1+1}). 
\endgather 
$$ 
Then  
$$ 
(\bar x_1,\dots,\bar x_{d-1};\bar y_1,\dots,\bar y_{d-1})= 
(\bar x(\mu^*);\bar y(\mu^*)). 
$$ 
 
We use the fact that  
$$ 
s_{\mu;a}(x(\mu);y(\mu))= 
s_{\mu;a}(\bar x(\mu);\bar y(\mu))= 
s_{\mu;a}(\bar x_1,\dots,\bar x_d;\bar y_1,\dots,\bar y_d) 
$$  
and employ the combinatorial formula \tht{4.5} and Theorem 4.6 to
compute the latter expression,  
$$ 
s_{\mu;a}(\bar x_1,\dots,\bar x_d;\bar y_1,\dots,\bar y_d)= 
\sum_\T\,\prod_{i=1}^d f_{\T^{-1}(i);a}(\bar x_i,\bar y_i),\tag5.8
$$ 
summed over all diagonal--strict tableaux of shape $\mu$. By
Proposition 4.3 and Theorem 4.6, this 
can also be written as follows 
$$ 
s_{\mu;a}(\bar x(\mu);\bar y(\mu))= 
\sum_{\nu\subset\mu} 
s_{\nu;a}(\bar x(\mu^*);\bar y(\mu^*)) 
f_{\mu/\nu;a}(\bar x_d,\bar y_d), 
$$ 
summed over all subdiagrams $\nu$ such that $\mu/\nu$ contains no 
$2\times2$ block of squares.  
 
By Theorem 5.1, $s_{\nu;a}(\bar x(\mu^*);\bar y(\mu^*))=0$ unless 
$\nu$ is contained in $\mu^*$. By the minimality property of $\mu^*$ 
mentioned above, this means that $\nu=\mu^*$. Repeating this argument we 
conclude that in the sum \tht{5.8}, there is only one tableau $\T=\T_0$ with 
(possibly) nonzero contribution: the chain of diagrams corresponding to 
$\T_0$ looks as follows: 
$$ 
\dots\subset(\mu^*)^*\subset\mu^*\subset\mu. 
$$ 
 
Finally, let us check that the contribution of $\T_0$ is indeed given 
by formula \tht{5.6}. Let $\mu[1]\subset\mu$ stand for the first diagonal 
hook in $\mu$, 
$$ 
\mu[1]=\{(1,j)\in\mu\}\cup\{(i,1)\in\mu\}. 
$$ 
We will prove that   
$$ 
f_{\mu/\mu^*;a}(\bar x_d,\bar y_d)= 
\prod_{(i,j)\in\mu[1]}(a'_{\mu_i-i+1/2}-a'_{-\mu'_j+j-1/2}). \tag5.9
$$ 
Moreover, the same argument will prove that  
$$ 
f_{\mu^*/(\mu^*)^*;a}(\bar x_{d-1},\bar y_{d-1})= 
\prod_{(i,j)\in\mu[2]}(a'_{\mu_i-i+1/2}-a'_{-\mu'_j+j-1/2}), 
$$ 
where $\mu[2]$ stands for the second diagonal hook, and so on. This 
will imply \tht{5.6}. 
 
By definition \tht{4.4}, 
$$ 
f_{\mu/\mu^*;a}(\bar x_d,\bar y_d)= 
(u+v) 
\prod\Sb k,\\ \text{$s_k$ is vertical}\endSb(u-a'_{\de_k-\ep_k}) 
\prod\Sb k,\\ \text{$s_k$ is horizontal}\endSb(v+a'_{\de_k-\ep_k}),\tag5.10
$$ 
where $(\ep_k,\de_k)$ are the midpoints of the interior sides 
$s_k$ of the shape $\mu/\mu^*$ and 
$$ 
u=\bar x_d=a_{p_1+1}=a'_{p_1+1/2}, \quad 
v=\bar y_d=\wha_{q_1+1}=-a'_{-(q_1+1/2)}\,. 
$$ 
 
We establish a bijective correspondence $s\leftrightarrow(i,j)$ 
between the sides $s=s_k$ and the squares $(i,j)\in\mu[1]$, except 
the diagonal square $(1,1)$, as follows: 
$$ 
s=(\ep,\de)\leftrightarrow(i,j)= 
\cases (1,\de+1), & \text{if $s$ is vertical,}\\ 
(\ep+1,1), & \text{if $s$ is horizontal}. \endcases 
$$ 
Note that if $(i,j)=(1,\de+1)$ then $\ep+1/2=\mu'_j$, and if 
$(i,j)=(\ep+1,1)$ then $\de+1/2=\mu_i$.  It follows that, under the above 
correspondence, the contribution of $s$ to \tht{5.10} coincides  
with the contribution of $(i,j)$ to  \tht{5.9}. 
As for the factor $(u+v)$ in \tht{5.10}, it 
coincides with the contribution of the square $(1,1)$. This proves 
\tht{5.9} and concludes the proof of the theorem. \qed  
\enddemo 
 
\proclaim{Corollary 5.3} If the numbers $a_i$ are pairwise 
distinct then $s_{\mu;a}(x(\mu);y(\mu))\ne0$. 
\endproclaim 
 
\demo{Proof} Recall a well--known claim: the sets 
$\{\mu_i-i+1\}_{i=1}^\infty$ and $\{j-\mu'_j\}_{j=1}^\infty$ do not 
intersect (and, moreover, their union is the whole $\Z$), see, e.g.,
\cite{Ma1, I, (1.7)}. It follows that, under the assumption on $a$,
all the factors in the product \tht{5.6} do not vanish. \qed 
\enddemo 
 
In the special case $a_i=i-\frac12$, Theorem 5.2 means that for any 
Young diagram $\mu=(p_1,\dots,p_d\,|\,q_1,\dots,q_d)$, 
$$ 
Fs_\mu(p_1+\tfrac12,\dots,p_d+\tfrac12;q_1+\tfrac12,\dots,q_d+\tfrac12)= 
\prod_{(i,j)\in\mu}(\mu_i-i+\mu'_j-j+1),\tag5.11
$$ 
which is equal to 
$$ 
\prod_{(i,j)\in\mu}h(i,j)=\frac{|\mu|!}{\dim\mu}, 
$$ 
the product of the hook lengths. There exist other useful expressions
for this product, in particular, 
$$ 
\prod_{(i,j)\in\mu}h(i,j)=\frac{\prod\limits_{i=1}^d p_i!q_i! 
\prod\limits_{i,j=1}^d(p_i+q_j+1)} 
{\prod\limits_{1\le i<k\le d}(p_i-p_k)(q_i-q_k)} 
$$ 
(see, e.g., \cite{BR, 7.14.1}, where 
one must take $k=l=d$, or \cite{Ol}). The next proposition 
provides a generalization of the last identity.  
 
\proclaim{Proposition 5.4} For any Young diagram 
$\mu=(p_1,\dots,p_d\,|\,q_1,\dots,q_d)$, the product  
$$
\multline
\prod_{(i,j)\in\mu}(a_{\mu_i-i+1}-a_{j-\mu'_j})\\
=\prod_{(i,j)\in\mu}(a'_{\mu_i-i+1/2}-a'_{-\mu'_j+j-1/2})= 
\prod_{(i,j)\in\mu}(a_{\mu_i-i+1}+\wha_{\mu'_j-j+1}),
\endmultline \tag5.12
$$ 
which gives the value of $s_{\mu;a}(x(\mu);y(\mu))$, is equal to  
$$ 
\frac{\prod\limits_{i=1}^d(a_{p_i+1}-a_1)\dots(a_{p_i+1}-a_{p_i}) 
(\wha_{q_i+1}-\wha_1)\dots(\wha_{q_i+1}-\wha_{q_i}) 
\prod\limits_{i,j=1}^d(a_{p_i+1}+\wha_{q_j+1})} 
{\prod\limits_{1\le i<k\le d}(a_{p_i+1}-a_{p_k+1}) 
(\wha_{q_i+1}-\wha_{q_k+1})}\,.\tag$5.12'$
$$ 
\endproclaim 
 
Note that in the special case $a_i=i-\frac12$, this coincides with the 
right--hand side of \tht{5.11}.  
 
\demo{Proof} The part of the product \tht{5.12} corresponding to the 
squares $(i,j)\in\mu$ with $i,j\le d$ exactly coincides with the 
second product in the numerator of \tht{$5.12'$}. It follows that 
the equality of 
\tht{5.12} and \tht{$5.12'$} can be reduced to the following identity
$$ 
\gather 
\prod\Sb(i,j)\in\mu \\ i\le d, j>d \endSb 
(a_{\mu_i-i+1}-a_{j-\mu'_j}) 
\cdot \prod\Sb(i,j)\in\mu \\ j\le d, i>d \endSb 
(\wha_{\mu'_j-j+1}-\wha_{i-\mu_i})\\ 
= 
\frac{\prod\limits_{i=1}^d 
(a_{p_i+1}-a_1)\dots(a_{p_i+1}-a_{p_i})} 
{\prod\limits_{1\le i<k\le d} 
(a_{p_i+1}-a_{p_k+1})}\cdot 
\frac{\prod\limits_{i=1}^d 
(\wha_{q_i+1}-\wha_1)\dots(\wha_{q_i+1}-\wha_{q_i})} 
{\prod\limits_{1\le i<k\le d} 
(\wha_{q_i+1}-\wha_{q_k+1})}\,. 
\endgather 
$$ 
 
By symmetry, it suffices to prove the identity 
$$ 
\prod\Sb(i,j)\in\mu \\ i\le d, j>d \endSb 
(a_{\mu_i-i+1}-a_{j-\mu'_j})= 
\frac{\prod\limits_{i=1}^d 
(a_{p_i+1}-a_1)\dots(a_{p_i+1}-a_{p_i})} 
{\prod\limits_{1\le i<k\le d} 
(a_{p_i+1}-a_{p_k+1})}\,,\tag5.13
$$ 
which reduces to the following claim: for any diagram 
$\mu$ and any $i=1,\dots,d$, 
$$ 
\{1,\dots,\mu_i-i\}=\{\mu_j-j+1\}_{j=i+1}^d\, \sqcup 
\{j-\mu'_j\}_{j=d+1}^{\mu_i}\,,\tag5.14
$$ 
a disjoint union of two sets.  
 
By virtue of the combinatorial fact mentioned in the proof of 
Corollary 5.2, the two sets in the right--hand side of \tht{5.14} are indeed 
disjoint.
 
Next, each of the two sets in the right--hand side of \tht{5.14} 
is contained in the set from the left--hand side. Indeed, 
the minimal element of the first set is $\mu_d-d+1\ge1$ and the 
maximal element is $\mu_{i+1}-i\le\mu_i-i$. Similarly, in the second 
set, the minimal element is $d+1-\mu'_{d+1}\ge d+1-d=1$ and the 
maximal element is $\mu_i-\mu'_{\mu_i}\le \mu_i-i$. Finally, the 
total number of elements in both sets is 
$(d-i)+(\mu_i-d)=\mu_i-i$, which is equal to the number of elements in 
the left--hand side. This proves \tht{5.14} and concludes the proof of the 
proposition. \qed 
\enddemo 
 
The next two results are similar to the characterization theorems for shifted 
Schur functions, cf. \cite{Ok1}, \cite{OO1, \S3}. 
 
\proclaim{Theorem 5.5 {\rm(Characterization theorem I)}} Let $\mu$ be 
an arbitrary Young diagram and  
$$  
D(\mu)=\{\la\in\Y \bigm| |\la|\le|\mu|, \la\ne\mu\}. 
$$ 
Assume that the numbers $a_i$ are pairwise distinct. Then, as an 
element of $\La$, $s_{\mu;a}$ is uniquely determined by the  
following two properties{\rm:} first, its top degree homogeneous component 
is the Schur function $s_\mu${\rm;} second, it vanishes at 
$(x(\la),y(\la))$ for all $\la\in D(\mu)$. 
\endproclaim 
 
\demo{Proof} Assume $F$ is a symmetric function with the same two 
properties. Then there exists an expansion of the form 
$$ 
F=s_{\mu;a} + \sum_{\nu\in D(\mu)}c(\nu) s_{\nu;a} 
$$ 
with certain numerical coefficients $c(\nu)$. We must prove that these 
coefficients are actually equal to zero. Let $X$ stand for the set of 
$\nu$'s with $c(\nu)\ne0$. Assume that $X$ is nonempty and choose a minimal 
diagram $\la\in X$ with respect to the partial ordering by inclusion. 
We get $F(x(\la);y(\la))=s_{\mu;a}(x(\la);y(\la))=0$, because $\la\in 
D(\mu)$. On the other hand, since $\la$ is minimal, we also have 
$s_{\nu;a}(x(\la);y(\la))=0$ for any $\nu\in X\setminus\{\la\}$, 
because $\nu$ is not contained in $\la$. But 
$s_{\la;a}(x(\la);y(\la))\ne0$, because of Corollary 5.3. This leads 
to a contradiction. \qed 
\enddemo 
 
\proclaim{Theorem 5.6 {\rm(Characterization theorem II)}} Let $\mu$ be 
an arbitrary Young diagram and $D(\mu)$ be as above. 
Assume again that the numbers $a_i$ are pairwise distinct. Then, as an 
element of $\La$, $s_{\mu;a}$ is uniquely determined by the  
following three properties{\rm:} first, its degree is less than or equal to 
$|\mu|${\rm;} second, its value at $(x(\mu);y(\mu))$ is given by formula
\tht{5.6}{\rm;} third, it vanishes at $(x(\la),y(\la))$ for all $\la\in D(\mu)$. 
\endproclaim 
 
\demo{Proof} Let $F\in\La$ possess the same three properties. Then  
$$ 
F=\sum_{\nu\in D(\mu)\cup\{\mu\}}c(\nu) s_{\nu;a} 
$$ 
with certain numerical coefficients $c(\nu)$. The same argument as 
above proves that $c(\nu)=0$ for all $\nu\in D(\mu)$, so that $F$ is 
proportional to $s_{\mu;a}$. Then the second property implies that 
$F$ is exactly equal to $s_{\mu;a}$. \qed 
\enddemo 
 
In the particular case of the Frobenius--Schur functions, the results
of this section are equivalent to similar claims for the shifted
Schur functions contained in \cite{Ok1}, \cite{OO1, \S3}.
Interpolation of arbitrary polynomials in terms of Schubert
polynomials is discussed in \cite{L2}. A very general scheme of
Newton interpolation for symmetric polynomials is developed in
\cite{Ok2}. In the latter paper one can also find the references to earlier
works by Knop, Sahi, and Okounkov. 
 
\head \S6. Sergeev--Pragacz formula \endhead 
 
Recall first the conventional Sergeev--Pragacz formula. Let $m,n$ be 
arbitrary nonnegative integers and let $\mu$ be an arbitrary diagram not
containing the square $(m+1,n+1)$. Then 
$$ 
s_\mu(x_1,\dots,x_m;y_1,\dots,y_n)= 
\frac{\sum_{w\in\frak S_m\times\frak S_n}\ep(w)\,  
w[f_\mu(x_1,\dots,x_m;y_1,\dots,y_n)]} 
{V(x_1,\dots,x_m)V(y_1,\dots,y_n)}\,,  \tag6.1
$$ 
where 
$$ 
f_\mu(x_1,\dots,x_m;y_1,\dots,y_n)= 
\prod_{i=1}^m x_i^{(\mu_i-i)_++m-i}\cdot 
\prod_{j=1}^n y_j^{(\mu'_j-j)_++n-j}\cdot 
\prod\Sb i\le m,\, j\le n\\ (i,j)\in\mu\endSb (x_i+y_j). \tag6.2
$$ 
Here $\frak S_m\times\frak S_n$ is the product of two symmetric groups 
acting on polynomials in $x_1,\dots,x_m$ and $y_1,\dots,y_n$ by permuting 
separately the $x$'s and $y$'s; for  
$w=(w_1,w_2)\in\frak S_m\times\frak S_n$, the symbol $\ep(w)$ 
means $\sgn(w_1)\sgn(w_2)$; 
$$ 
V(x_1,\dots,x_m)=\prod_{1\le i<k\le m}(x_i-x_k), \qquad 
V(y_1,\dots,y_n)=\prod_{1\le j<l\le n}(y_j-y_l), 
$$ 
and  
$$ 
(k)_+=\max(k,0). 
$$ 
 
Note that for $n=0$ the formula \tht{6.1} reduces to the classical 
formula \tht{0.6} for the Schur polynomial $s_\mu(x_1,\dots,x_m)$. Proofs of 
\tht{6.1}--\tht{6.2} can be found in \cite{PT}, \cite{Ma1, I.3, Ex. 23}. 
 
Our aim in this section is to establish an analog of the formulas 
\tht{6.1}--\tht{6.2} for the multiparameter supersymmetric Schur polynomials 
depending on a sequence of parameters $a=(a_i)$. According to our
basic principle, we will deal with equal number of the $x$'s and $y$'s. 
 
\proclaim{Theorem 6.1} Let $n=1,2,\dots$ and let $\mu$ be an arbitrary 
Young diagram such that $d=d(\mu)\le n$. Then, in the above notation, 
$$ 
s_{\mu;a}(x_1,\dots,x_n;y_1,\dots,y_n)= 
\frac{\sum_{w\in\frak S_n\times\frak S_n}\ep(w)\,  
w[f_{\mu;a}(x_1,\dots,x_n;y_1,\dots,y_n)]} 
{V(x_1,\dots,x_n)V(y_1,\dots,y_n)}\,,  \tag6.3
$$ 
where 
$$ 
\gather 
f_{\mu;a}(x_1,\dots,x_n;y_1,\dots,y_n)= 
\prod_{i=1}^d (x_i\,|\,a)^{\mu_i-i}x_i^{(n-\mu_i)_+} 
(y_i\,|\,\wha)^{\mu'_i-i}y_i^{(n-\mu'_i)_+}\\ 
\times\prod_{i=d+1}^n x_i^{n-i} y_i^{n-i}\cdot 
\prod\Sb i,j\le n\\ (i,j)\in\mu\endSb (x_i+y_j) \tag6.4
\endgather 
$$ 
and 
$$ 
(x\,|\,a)^m=\cases (x-a_1)\dots(x-a_m), & m\ge1, \\ 
0, & m=0. \endcases 
$$ 
\endproclaim 
 
Note that for the sequence $a\equiv0$ formulas \tht{6.3}--\tht{6.4} reduce to 
formulas \tht{6.1}--\tht{6.2} with $m=n$. Indeed, it suffices to check that 
$f_{\mu;a}$ reduces to $f_\mu$. To do this, let us compare \tht{6.2} 
and \tht{6.4}. The last product in both formulas is the same. Then 
we remark that  
$$ 
(\mu_i-n)_+ +n-i=\cases \mu_i-i+(n-\mu_i)_+, & i=1,\dots, d,\\ 
n-i, &i=d+1,\dots, n, \endcases 
$$ 
and likewise for $\mu'$, which implies the equality 
$f_\mu=f_{\mu;a}\bigm|_{a\equiv0}$.  
 
\demo{Proof} The proof of Theorem 6.1 is divided into three lemmas. 
 
\proclaim{Lemma 6.2} The expression \tht{6.3} is a supersymmetric 
polynomial.  
\endproclaim 
 
\demo{Proof} One may argue exactly as in the proof of
Proposition 2.3 in \cite{PT} (see also \cite{Ma1, I.3, Ex. 24}). For the
reader's convenience we present the argument. 
 
Obviously, \tht{6.3} is a polynomial, which is separately 
symmetric in the $x$'s and $y$'s. Let us verify the cancellation 
property: assume that $x_i=t=-y_j$ for certain indices 
$i,j$, and let us prove that \tht{6.3} does not depend on $t$. Consider the 
expression \tht{6.3} as a function in $t$. This is a rational function, which is 
actually a polynomial. The degree of the denominator 
is exactly $2n-2$, so that it suffices to prove that the numerator 
has degree less than or equal to $2n-2$. Next, it suffices to prove the 
latter claim for the polynomial \tht{6.4}. 
 
Without loss of generality, one may assume $i\le j$. Consider three 
cases.  
 
First case: $i\le j\le d$. Then $(i,j)\in\mu$, so that \tht{6.4} 
contains the factor $(x_i+y_j)$. This factor vanishes when 
$x_i=-y_j$, consequently, \tht{6.4} vanishes identically.  
 
Second case: $i\le d<j\le n$. The degree of the expression \tht{6.4} 
with respect to $t$ is less than or equal to the sum of three terms: 
$$ 
(\mu_i-i+(n-\mu_i)_+)+(n-j)+(\mu_i+\mu'_j). 
$$ 
Again, if $(i,j)\in\mu$ then \tht{6.4} is identically equal to 0, so 
that one may assume $(i,j)\notin\mu$. This means that $\mu_i\le j-1$ and 
$\mu'_j\le i-1$; in particular, $\mu_i<n$. It follows that the above sum
reduces to  
$$ 
(n-i)+(n-j)+(\mu_i+\mu'_j), 
$$ 
which is $\le 2n-2$, because $\mu_i+\mu'_j\le i+j-2$, as was 
mentioned above.  
 
Third case: $d+1\le i\le j\le n$. Then the degree in question is 
$$ 
(n-i)+(n-j)+(\mu_i+\mu'_j)\le 2n-i-j+(j-1)+(i-1)=2n-2. 
$$ 
This concludes the proof of the lemma. \qed 
\enddemo 
 
Let us temporarily denote the right--hand side of \tht{6.3} by 
$s'_{\mu;a}(x_1,\dots,x_n;y_1,\dots,y_n)$. 
 
\proclaim{Lemma 6.3} We have  
$$
\multline 
s'_{\mu;a}(x_1,\dots,x_n;y_1,\dots,y_n)\bigm|_{x_n=y_n=0}\\
=\cases 0, & d=n,\\ 
s'_{\mu;a}(x_1,\dots,x_{n-1};y_1,\dots,y_{n-1}), & d\le n-1. \endcases 
\endmultline
$$ 
\endproclaim 
 
\demo{Proof} a) If $d=n$ then the polynomial 
$f_{\mu;a}(x_1,\dots,x_n;y_1,\dots,y_n)$ is divisible by 
$$ 
\prod_{i,j=1}^n(x_i+y_j), 
$$ 
and the same holds for its transformations by elements $w$. Since the 
above product vanishes when $x_n=y_n=0$, the whole expression 
$s'_{\mu;a}(x_1,\dots,x_n;y_1,\dots,y_n)$ vanishes, too.  
 
b) Next, let us assume that $d\le n-1$ and let us prove the following 
claim:  
$$ 
f_{\mu;a}(x_1,\dots,x_n;y_1,\dots,y_n)\bigm|_{x_i=y_j=0}=0 
\qquad \text{unless $i=j=n$.} 
$$ 
 
Indeed, vanishing holds if $(i,j)\in\mu$, because then the polynomial in 
question is divisible by $x_i+y_j$. Assume $i\le d$. If 
$\mu_i\ge n$ then $(i,j)\in\mu$, which implies vanishing. If 
$\mu_i<n$ then $(n-\mu_i)_+>0$, so that vanishing holds thanks to the 
factor $x_i^{(n-\mu_i)_+}$ in the expression \tht{6.4}. Similar 
argument also holds when $j\le d$. 
 
Thus, we may assume that both $i> d$ and $j> d$. Then the 
expression \tht{6.4} contains the factor $x_i^{n-i}y_j^{n-j}$, so that 
the only possibility of nonvanishing may occur for $i=j=n$.  
 
c) The numerator of \tht{6.3} is an alternate sum of terms indexed by couples of 
permutations $w=(w_1,w_2)$. If $w$ does not fix the 
indeterminates $(x_n,y_n)$, which are specialized to zero, then the 
corresponding term vanishes after the specialization because of the 
claim b) proved above. Consequently, only terms indexed by elements 
$w$ fixing $(x_n,y_n)$ may give a nonzero contribution. Note that 
these are actually elements of the group  
$\frak S_{n-1}\times\frak S_{n-1}$.  
 
On the other hand, note that 
$$ 
\gather 
V(x_1,\dots,x_{n-1},0)V(y_1,\dots,y_{n-1},0)= 
V(x_1,\dots,x_{n-1})V(y_1,\dots,y_{n-1})\\ 
\times\, x_1\dots x_{n-1}y_1\dots y_{n-1}\,. 
\endgather 
$$ 
 
d) Thus, it remains to check that  
$$ 
\gather 
f_{\mu;a}(x_1,\dots,x_{n-1},0;y_1,\dots,y_{n-1},0) 
=f_{\mu;a}(x_1,\dots,x_{n-1};y_1,\dots,y_{n-1})\\ 
\times\, x_1\dots x_{n-1}y_1\dots y_{n-1}\,. 
\endgather 
$$ 
Let us examine the behavior of the expression \tht{6.4} under 
the specialization $x_n=y_n=0$. The right--hand side of \tht{6.4} consists of 
three 
products. The second product, which is equal to  
$$ 
\prod_{i=d+1}^n x_i^{n-i} y_i^{n-i}\,, 
$$ 
turns into the similar expression for $n-1$  
multiplied by 
$$ 
\prod_{i=d+1}^{n-1} x_i y_i. 
$$ 
Consequently, we must prove that the remaining expression in \tht{6.4}, 
which is equal to  
$$ 
\prod_{i=1}^d (x_i\,|\,a)^{\mu_i-i}x_i^{(n-\mu_i)_+} 
(y_i\,|\,\wha)^{\mu'_i-i}y_i^{(n-\mu'_i)_+} 
\cdot  
\prod\Sb i, j\le n\\ (i,j)\in\mu\endSb (x_i+y_j), 
$$ 
turns into the similar expression for $n-1$ multiplied by  
$$ 
\prod_{i=1}^d x_iy_i. 
$$ 
 
This is equivalent to the following claim: for any $i=1,\dots,d$, 
$$ 
(n-\mu_i)_++\ep_i=(n-1-\mu_i)_++1, \tag6.5
$$ 
where  
$$ 
\ep_i=\cases 1, & (i,n)\in\mu,\\ 
0, & (i,n)\notin\mu, \endcases 
$$ 
and similarly for $\mu'$.  
 
Consider two cases: $\mu_i\ge n$ and $\mu_i\le n-1$. In the former 
case, $(n-\mu_i)_+=(n-1-\mu_i)_+=0$ and $\ep_i=1$, so that \tht{6.5} 
holds. In the latter case, $(n-\mu_i)_+=(n-1-\mu_i)_++1$ and 
$\ep_i=0$, so that \tht{6.5} is again true. This concludes the proof. 
\qed 
\enddemo 
 
By Lemma 6.2 and Lemma 6.3, there 
exists an element $s'_{\mu;a}\in\La$, of degree $|\mu|$ and  such that the 
corresponding 
supersymmetric polynomials in $n+n$ variables are given by the 
expressions \tht{6.3} when $n\ge d$ and vanish when $n<d$. We aim at
proving that $s'_{\mu;a}=s_{\mu;a}$ by making use of the 
characterization theorems from \S5.
 
\proclaim{Lemma 6.4} The element $s'_{\mu;a}\in\La$ defined above
satisfies the vanishing condition of Theorem\/ {\rm5.5}, i.e., 
$$ 
s'_{\mu;a}(x(\la);y(\la))=0 \qquad  
\text{unless $\mu\subseteq\la$.} 
$$ 
\endproclaim 
 
\demo{Proof} We will see that the vanishing in question is ensured by 
the product 
$$ 
\prod_{i=1}^d (x_i\,|\,a)^{\mu_i-i}(y_i\,|\,\wha)^{\mu'_i-i} 
$$ 
entering the expression \tht{6.4}.  
 
Write both diagrams in the Frobenius notation,  
$$ 
\mu=(p_1,\dots,p_d\,|\,q_1,\dots,q_d), \qquad 
\la=(P_1,\dots,P_D\,|\,Q_1,\dots,Q_D), 
$$ 
and recall that 
$$ 
\gather 
x(\la)_i=a_{P_i+1},\quad  
y(\la)_i=\wha_{Q_i+1},\quad 1\le i\le D,\\  
x(\la)_i=y(\la)_i=0, \quad i>D. 
\endgather 
$$ 
By the definition of $s'_{\mu;a}$, vanishing holds if $D<d$, because 
$x(\la)$ and $y(\la)$ have at most $D$ nonzero coordinates, so that 
we may assume $D\ge d$.  
 
Abbreviate   
$$ 
x_i=x(\la)_i=a_{P_i+1},\quad  
y_i=y(\la)_i=\wha_{Q_i+1},\quad 1\le i\le D. 
$$ 
It suffices to prove that for any  
$w=(w_1,w_2)\in\frak S_D\times\frak S_D$, 
$$ 
w\left[\prod_{i=1}^d (x_i\,|\,a)^{p_i}(y_i\,|\,\wha)^{q_i}\right]=0 
$$ 
unless $P_i\ge p_i$ and $Q_i\ge q_i$ for all $i=1,\dots,d$. By the 
symmetry $x\leftrightarrow y$, $a\leftrightarrow \wha$, it suffices 
to prove that for any $w_1\in\frak S_D$, 
$$ 
w_1\left[\prod_{i=1}^d (x_i\,|\,a)^{p_i}\right]=0 \tag6.6
$$ 
unless $P_i\ge p_i$ for all $i=1,\dots,d$. 
 
The left--hand side in \tht{6.6} has the form 
$$ 
\prod_{i=1}^d (a_{P_{j_i}+1}\,|\,a)^{p_i}= 
\prod_{i=1}^d (a_{P_{j_i}+1}-a_1)\dots(a_{P_{j_i}+1}-a_{p_i})\,, 
$$ 
where $j_1,\dots,j_d$ is a certain $d$--tuple of pairwise distinct 
indices from $\{1,\dots,D\}$. If \tht{6.6} does not hold then 
$$
P_{j_1}\ge p_1, \dots, P_{j_d}\ge p_d.\tag6.7
$$
Let us check that this implies 
$$
P_1\ge p_1,\dots, P_d\ge p_d. \tag 6.8
$$
Indeed, recall that 
$$
P_1>\dots>P_d, \qquad p_1>\dots>p_d.
$$
Together with \tht{6.7} this implies that among the numbers
$P_1,\dots,P_d$, there is at least one number $\ge p_1$ (namely,
$P_{j_1}$), at least two numbers $\ge p_2$ (namely, $P_{j_1},
P_{j_2}$), and so on, which implies \tht{6.8}.  \qed 
\enddemo 
 
Now we are in a position to prove the equality 
$s'_{\mu;a}=s_{\mu;a}$. Comparing the formula \tht{6.3}--\tht{6.4} defining 
$s'_{\mu;a}$ with the Sergeev--Pragacz formula \tht{6.1}--\tht{6.2} for the 
Schur 
function $s_\mu$ we see that the top degree homogeneous component of 
$s'_{\mu;a}$ coincides with $s_\mu$. By Lemma 6.4, $s'_{\mu;a}$ 
possesses the same vanishing property as $s_{\mu;a}$. Consequently, 
by Theorem 5.5, $s'_{\mu;a}=s_{\mu;a}$. Note that in Theorem 5.5, the 
numbers $a_i$ are required to be pairwise distinct but one may attain 
this by making use of the continuity argument, because both 
$s'_{\mu;a}$ and $s_{\mu;a}$ depends on the parameters continuously. 
 
An alternative way is to apply Theorem 5.6. Then we must verify for 
$s'_{\mu;a}$ the three properties listed in the statement of this 
theorem. The first property (control of degree) is obvious. The third 
property (vanishing) is ensured by Lemma 6.4. The second property 
(required value at $(x(\mu);y(\mu))$) is verified as follows. {}From 
the proof of Lemma 6.4 one sees that in formula \tht{6.3} applied to the 
variables $(x(\mu);y(\mu))$, all terms with $w\ne e$ are zero. Next, 
it is readily seen that the term with $w=e$ leads to the expression 
\tht{$5.12'$}, which, by Proposition 5.4, coincides with 
$s_{\mu;a}(x(\mu);y(\mu))$.   
 
This concludes the proof of Theorem 6.1. \qed 
\enddemo 
 
The following result is an analog of the Berele--Regev 
factorization property for the supersymmetric Schur functions 
\cite{BR}, \cite{Ma1, I.3, Ex. 23}.  
 
\proclaim{Corollary 6.5} Let $\mu=(p_1,\dots,p_d\,|\,q_1,\dots,q_d)$ 
be a Young diagram of depth $d$. Then 
$$ 
s_{\mu;a}(x_1,\dots,x_d;y_1,\dots,y_d)= 
\frac{\det[(x_i\,|\,a)^{p_j}]_{i,j=1}^d}{V(x_1,\dots,x_d)}\cdot 
\frac{\det[(y_i\,|\,\wha)^{q_j}]_{i,j=1}^d}{V(y_1,\dots,y_d)}\cdot 
\prod_{i,j=1}^d (x_i+y_j)\,. 
$$ 
\endproclaim 
 
\demo{Proof} We apply formulas \tht{6.3} and \tht{6.4} for $n=d$. 
Then the last  
product in \tht{6.4} coincides with the product above, which is 
invariant under permutations $w$.  Consequently, the alternate sum in 
\tht{6.3} becomes the product of two determinants. Note also that for 
each $i=1,\dots,d$, we have $\mu_i-i=p_i$,  
$(n-\mu_i)_+=(d-\mu_i)_+=0$, and, similarly, $\mu'_i-i=q_i$, 
$(n-\mu'_i)_+=0$. \qed 
\enddemo

\head 7. Transition coefficients \endhead 
 
Formulas \tht{3.1} and \tht{4.9} yield the generating series for the one--row 
and one--column multiparameter Schur functions; let us rewrite them in 
slightly different notation: 
$$ 
\gathered
1+\sum_{p=0}^\infty\frac{s_{(p\,|\,0);a}}{(u\,|\,a)^{p+1}}=H(u)= 
1+\sum_{p=0}^\infty\frac{s_{(p\,|\,0)}}{u^{p+1}}\,, \\ 
1+\sum_{q=0}^\infty\frac{s_{(0\,|\,q);a}}{(v\,|\,\wha)^{q+1}}=E(v)= 
1+\sum_{q=0}^\infty\frac{s_{(0\,|\,q)}}{v^{q+1}}\,. 
\endgathered \tag7.1 
$$ 
 
The next proposition yields the generating series for the hook functions 
$s_{(p\,|\,q);a}$; this series is an element of $\La[[u^{-1},v^{-1}]]$.  
 
\proclaim{Proposition 7.1} We have 
$$ 
1+(u+v)\sum_{p,q=0}^\infty 
\frac{s_{(p\,|\,q);a}} 
{(u\,|\,a)^{p+1}(v\,|\,\wha)^{q+1}} 
=H(u)E(v).  
$$ 
\endproclaim 
 
\demo{Proof} This is a generalization of Theorem 2.3, and we will argue as in 
the proof of that theorem. By virtue of \tht{7.1}, 
the equality in question is equivalent to  
$$ 
1+(u+v)\sum_{p,q=0}^\infty 
\frac{s_{(p\,|\,q);a}} 
{(u\,|\,a)^{p+1}(v\,|\,\wha)^{q+1}}= 
\left(1+\sum_{p=0}^\infty\frac{s_{(p\,|\,0);a}}{(u\,|\,a)^{p+1}}\right) 
\left(1+\sum_{q=0}^\infty\frac{s_{(0\,|\,q);a}}
{(v\,|\,\wha)^{q+1}}\right). 
$$ 
Using the identity 
$$ 
\gather 
\frac{u+v}{(u\,|\,a)^{p+1}(v\,|\,\wha)^{q+1}}= 
\frac{(u-a_{p+1})+(v-\wha_{q+1})+(a_{p+1}+\wha_{q+1})} 
{(u\,|\,a)^{p+1}(v\,|\,\wha)^{q+1}}\\ 
=\frac1{(u\,|\,a)^{p}(v\,|\,\wha)^{q+1}}+ 
\frac1{(u\,|\,a)^{p+1}(v\,|\,\wha)^{q}}+ 
\frac{a_{p+1}+\wha_{q+1}}{(u\,|\,a)^{p+1}(v\,|\,\wha)^{q+1}} 
\endgather 
$$ 
we reduce this to the following system of relations, 
where $p,q=0,1,\dots$: 
$$ 
s_{(p+1\,|\,q);a}+s_{(p\,|\,q+1);a}+ 
(a_{p+1}+\wha_{q+1})s_{(p\,|\,q);a}=s_{(p\,|\,0);a}s_{(0\,|\,q);a}\,. 
$$ 
 
By virtue of \tht{1.4}, the expansion of the product 
$s_{(p\,|\,0);a}s_{(0\,|\,q);a}$ into 
a linear combination of the functions $s_{\nu;a}$ has the form 
$$ 
s_{(p\,|\,0);a}s_{(0\,|\,q);a}= 
s_{(p+1\,|\,q);a}+s_{(p\,|\,q+1);a}+ 
\sum_{\nu:\, |\nu|\le p+q+1}c(\nu)s_{\nu;a}\,. 
$$ 
Let $X$ be the set of those $\nu$'s which enter this sum with 
nonzero coefficients $c(\nu)$. We claim that $X$ contains at most the 
diagram $(p\,|\,q)$; here we will use the fact that 
$|\nu|\le p+q+1$. Indeed, let $\la$ be a minimal (with respect to 
inclusion) diagram in $X$ and let us evaluate both sides at 
$(x(\la);y(\la))$. On the right, the result is nonzero, because 
$s_{\la;a}(x(\la);y(\la))\ne0$ \footnote{Here we tacitly assume that 
the sequence $a=(a_i)$ has no repetitions in order to apply Corollary 
5.3. To cover the case when repetitions are present, one can use the 
continuity argument.} while all other terms on the right have zero 
contributions, because neither a diagram $\nu\ne\la$ from $X$ nor 
$(p+1\,|\,q)$ and $(p\,|\,q+1)$ are contained in $\la$.  So, the result 
of the evaluation on the left is nonzero, too. This implies that $\la$ 
contains both $(p\,|\,0)$ and $(0\,|\,q)$, which is only possible for 
$\la=(p\,|\,q)$.  
 
Thus, the above expansion takes the form  
$$ 
s_{(p\,|\,0);a}s_{(0\,|\,q);a}= 
s_{(p+1\,|\,q);a}s_{(p\,|\,q+1);a}+ 
\operatorname{const}\,s_{(p\,|\,q);a}\,. 
$$ 
Evaluating both sides at 
$$ 
(x(p\,|\,q);y(p\,|\,q))=(a_{p+1},0,0,\dots;\wha_{q+1},0,0,\dots) 
$$ 
we get 
$$ 
\operatorname{const}=\frac{s_{(p\,|\,0);a}(x(p\,|\,q);y(p\,|\,q))\, 
s_{(0\,|\,q);a}(x(p\,|\,q);y(p\,|\,q))} 
{s_{(p\,|\,q)}(x(p\,|\,q);y(p\,|\,q))}\,. 
$$ 
 
The right--hand side can be readily evaluated by making use of Corollary 6.5. 
This 
gives the desired result: $\operatorname{const}=a_{p+1}+\wha_{q+1}$. 
\qed  
\enddemo

Let $b=(b_i)_{i\in\Z}$ be another sequence of parameters. We aim at
expressing the functions $s_{\mu;a}$ through the functions 
$s_{\nu;b}$. In particular, we are interested in the expansion on 
the conventional supersymmetric Schur functions $s_\nu$, which 
correspond to $b\equiv0$.  
 
For $p\ge p'\ge0$ set 
$$ 
c_{pp'}(a,b)=h_{p-p'}(b_1,\dots,b_{p'+1};-a_1,\dots,-a_p), \tag7.2 
$$ 
where $h_{p-p'}$ is the conventional supersymmetric $h$-function 
of degree $p-p'$. In particular, 
$$ 
c_{pp'}(a,0)=h_{p-p'}(0,\dots,0;-a_1,\dots,-a_p) 
=(-1)^{p-p'}e_{p-p'}(a_1,\dots,a_p), \tag7.3 
$$ 
cf. \tht{2.9}. Note that $c_{pp}(a,b)=1$. 
 
\proclaim{Proposition 7.2 {\rm (cf. Proposition 2.4)}} We have 
$$ 
s_{(p\,|\,q);a}=\sum_{p'=0}^p\sum_{q'=0}^q  
c_{pp'}(a,b)\,c_{qq'}(\wha,\whb) \,s_{(p'\,|\,q');b}. \tag7.4 
$$ 
\endproclaim 
 
\demo{Proof} We have the identity 
$$ 
(u+v)\sum_{p,q=0}^\infty 
\frac{s_{(p\,|\,q);a}} 
{(u\,|\,a)^{p+1}(v\,|\,\wha)^{q+1}} 
=(u+v)\sum_{p,q=0}^\infty 
\frac{s_{(p\,|\,q);b}} 
{(u\,|\,b)^{p+1}(v\,|\,\whb)^{q+1}}\,, 
$$  
because, by Proposition 7.1, both sides are equal to 
$H(u)E(v)-1$. Then we argue exactly as in the proof of Proposition 
2.4. \qed 
\enddemo 
 
\proclaim{Theorem 7.3 {\rm(cf. Theorem 2.6)}} Let $a=(a_i)_{i\in\Z}$ and 
$b=(b_i)_{i\in\Z}$ 
be two sequences of parameters. In the expansion 
$$ 
s_{\mu;a}=\sum_\nu c_{\mu\nu}(a,b)\,s_{\nu;b} \tag7.5 
$$ 
the coefficients $c_{\mu\nu}(a,b)$ vanish unless $\nu\subseteq\mu$ 
and $d(\nu)=d(\mu)$.  
 
Assume $d(\nu)=d(\mu)$ and write both diagrams in the Frobenius 
notation,  
$$ 
\mu=(p_1,\dots,p_d\,|\,q_1,\dots,q_d),\quad 
\nu=(p'_1,\dots,p'_d\,|\,q'_1,\dots,q'_d). 
$$ 
Then   
$$ 
c_{\mu\nu}(a,b)=\det[c_{p_i,p'_j}(a,b)]\, 
\det[c_{q_i,q'_j}(\wha,\whb)], \tag7.6 
$$ 
where the determinants are of order $d$ and the coefficients in the 
right--hand side are defined by \tht{7.2}.  
\endproclaim 
 
\demo{Proof} The same argument as in the proof of Theorem 2.6 reduces 
the desired claim to the special case $d=1$, which was the subject of 
Proposition 7.2. \qed  
\enddemo

\head Appendix by Vladimir Ivanov: Proof of the combinatorial formula 
for multiparameter Schur functions \endhead 

The aim of this appendix is to prove Theorem 4.6, which is restated
below as Theorem A.6. Here we are using a slightly different notation
for the combinatorial sum, the $h$--functions and the $s$--functions.
 
\definition{Definition A.1} 
Suppose $(a_i)_{i\in \Z}$ 
is an arbitrary sequence of complex numbers. 
For $k=0,1,2,\dots$, set 
$$ 
(x\mid a)^k=\prod_{i=1}^k(x-a_i).
$$ 
\enddefinition 
 
\definition{Definition A.2} 
Fix $n=1,2,\dots$ and denote by $\A_n$ the 
ordered alphabet $\{1'<1<2'<2<\dots<n'<n\}$. 
Put $|i'| =|i|=i$ for $i=1,2,\dots,n$. 
\enddefinition 
 
\definition{Definition A.3} 
Fix a skew Young diagram $\la/\mu$. 
A map $T:\la/\mu\to \A_n$ 
is called {\it a tableau of 
shape $\la/\mu$ and of order $n$} 
if the following conditions hold: 
\roster 
\item $T(i,j)\le T(i,j+1)$; 
\item $T(i,j)\le T(i+1,j)$; 
\item for each $i=1,2,\dots,n$, 
there is at most one symbol $i'$ 
in each row; 
\item for each $i=1,2,\dots,n$, 
there is at most one symbol $i$ 
in each column. 
\endroster 
Let us denote by $\Tab(\la/\mu,n)$ the set 
of all tableaux of shape $\la/\mu$ 
and order $n$. 
\enddefinition 
 
\definition{Definition A.4} 
Suppose that $x_1,x_2,\dots,y_1,y_2,\dots$ 
are variables, $(a_i)_{i\in \Z}$ 
is an arbitrary sequence, 
$\mu\subset\la$ are Young diagrams. Set 
$$ 
\multline
\Sigma_{\la/\mu\mid n}(x;y\mid a)\\
= \sum_{T\in\Tab(\la/\mu,n)} 
\prod\Sb \ssq\in\la/\mu \\ T(\ssq)=1,2,\dots,n 
\endSb (x_{T(\ssq)}-a_{c(\ssq)}) 
\prod\Sb \ssq\in\la/\mu \\ T(\ssq)=1',2',\dots,n' 
\endSb (y_{T(\ssq)}+a_{c(\ssq)}), 
\endmultline
$$ 
where $c(i,j)=j-i$, cf. (4.7). 
\enddefinition 
 
\definition{Definition A.5} 
Let $\mu\subset\la$. Set 
$$ 
\multline 
s_{\la/\mu}(x_1,x_2,\dots,x_n;y_1,y_2,\dots,y_n\mid a)=\\ 
\det[h_{\la_i-\mu_j+j-i}(x_1,x_2,\dots,x_n; 
y_1,y_2,\dots,y_n\mid \tau^{\mu_j-j+1}a)]_{1\le i,j\le \ell(\la)}, 
\endmultline 
$$ 
where the polynomials $h_k(x_1,\dots,x_n;y_1,\dots,y_n\mid a)$ are given
by the generating series
$$
1+\sum_{k=1}^\infty\frac{h_k(x_1,\dots,x_n;y_1,\dots,y_n\mid a)}
{(u\mid a)^k}=\prod_{i=1}^n\frac{u+y_i}{u-x_i}\,,
$$
$h_k\equiv 0$ if $k<0$, and $h_0\equiv 1$, cf. (3.4).
\enddefinition 
 
\proclaim{Theorem A.6} 
$$ 
s_{\la/\mu}(x_1,x_2,\dots,x_n;y_1,y_2,\dots,y_n\mid a)= 
\Sigma_{\la/\mu\mid n}(x;y\mid a). 
$$ 
\endproclaim 

\demo{Proof} 
First let us prove this claim
when $\la=(k),\mu=\varnothing$. 
Then $s_{\la/\mu}(x;y\mid a)= h_k(x;y\mid a)$ 
and our claim is equivalent 
to the following branching rule (see Proposition 4.3): 
$$ 
\multline 
h_k(x_1,\dots,x_{n-1},x_n;y_1,\dots,y_{n-1},y_n\mid a)= 
h_k(x_1,\dots,x_{n-1};y_1,\dots,y_{n-1}\mid a)\\ 
+(x_n+y_n)\sum_{r=0}^{k-1}\,h_r(x_1,\dots,x_{n-1};y_1,\dots,y_{n-1}\mid a)\,
\prod_{t=1}^{k-1-r}(x_n-a_{k-t}). 
\endmultline 
$$ 
 
Let us denote 
$h_k(x_1,\dots,x_{n-1},x_n;y_1,\dots,y_{n-1},y_n\mid a)$ 
by $h_{k\mid n}(x;y\mid a)$. 
{}From the definition of $h_k$ we obtain 
$$ 
1+\sum_{k=1}^\infty\frac{h_{k\mid n}(x;y\mid a)}{(u\mid a)^k}= 
\prod_{i=1}^n\frac{u+y_i}{u-x_i}= 
\frac{u+y_n}{u-x_n}+ 
\sum_{k=1}^\infty\frac{h_{k\mid n-1}(x;y\mid a)}{(u\mid a)^k}
\frac{u+y_n}{u-x_n}\,. 
\tag A.1 
$$ 
For any $k=1,2,\dots$ we have 
$$ 
\frac{u+y}{u-x}= 
1+(x+y)\sum_{r=1}^\infty 
\frac{(x\mid \tau^ka)^{r-1}}{(u\mid \tau^ka)^r}.\tag A.2 
$$ 
This follows, for example, from 
Molev's results \cite{Mo, Prop. 1.2, 
Theorem 2.1} if we put 
$m=n=1$. 
Using (A.1) and (A.2), we obtain 
$$ 
\multline 
1+\sum_{k=1}^\infty\frac{h_{k\mid n}(x;y\mid a)}{(u\mid a)^k}= 
\frac{u+y_n}{u-x_n}+ 
\sum_{k=1}^\infty\frac{h_{k\mid n-1}(x;y\mid a)}
{(u\mid a)^k}\frac{u+y_n}{u-x_n}\\ 
=1+(x_n+y_n)\sum_{k=1}^\infty\frac{(x_n\mid a)^{k-1}}{(u\mid a)^k}+ 
\sum_{k=1}^\infty\frac{h_{k\mid n-1}(x;y\mid a)}{(u\mid a)^k}\cdot 
\left(1+(x_n+y_n)\sum_{r=1}^\infty\frac{(x_n\mid \tau^ka)^{r-1}} 
{(u\mid \tau^ka)^r}\right)\\ 
=1+\sum_{k=1}^\infty 
\frac{h_{k\mid n-1}(x;y\mid a)+ 
(x_n+y_n)\sum_{r=0}^{k-1} 
h_{r\mid n-1}(x;y\mid a)\prod_{t=1}^{k-1-r}(x_n-a_{k-t})}
{(u\mid a)^k}. 
\endmultline 
$$ 
Therefore, when $\la=(k)$ and 
$\mu=\varnothing$, the assertion of the theorem is 
proved: 
$$ 
h_{k\mid n}(x;y\mid a)= 
\sum_{T\in\Tab((k),n)} 
\prod\Sb \ssq\in\la/\mu \\ T(\ssq)=1,2,\dots,n 
\endSb (x_{T(\ssq)}-a_{c(\ssq)}) 
\prod\Sb \ssq\in\la/\mu \\ T(\ssq)=1',2',\dots,n' 
\endSb (y_{T(\ssq)}+a_{c(\ssq)}). \tag A.3 
$$ 
 
Now suppose that $\la/\mu$ is an arbitrary skew Young diagram. 
Let us show that 
$$ 
\Sigma_{\la/\mu\mid n}(x;y\mid a)= 
\det[h_{\la_i-\mu_j+j-i\mid n}
(x;y\mid \tau^{\mu_j-j+1}a)]_{1\le i,j\le \ell(\la)}. 
$$ 
We use an appropriate modification of the Gessel-Viennot method (see
\cite{GV}, \cite{Sa, \S4.5} for a detailed description of this method). 

Suppose we have $2n$ variables  
$x_1,\dots,x_n,y_1,\dots,y_n$. 
We define a {\it path} 
as a finite sequence 
$p=(p(0),p(1),\dots, p(l))$ of points in the
$(x,y)$--plane, such that:

$\bullet$ Each point $p(t)$ belongs to $\Z\times\{0,1,2,\dots,n\}$.  

$\bullet$ The path starts at the horizontal line $y=0$ and ends on
the line $y=n$, i.e., $p(0)=(m_1,0)$ and $p(l)=(m_2,n)$ with some
$m_1, m_2\in\Z$.
  
$\bullet$ Each step $p(t+1)-p(t)$ is of the form $(1,0)$, $(1,1)$ or
$(0,1)$ (in accordance with these three options, 
we speak of a {\it horizontal, diagonal\/} or {\it vertical\/} step, 
respectively).

$\bullet$ The first step $p(1)-p(0)$ is either diagonal or vertical
but not horizontal.
 
We will deal with collections of paths 
$L=(p_i)_{1\le i\le \ell(\la)}$ 
such that $p_i$ starts 
at $(\mu_i-i,0)$, 
ends at $(\la_i-i,n)$, and 
$p_i\cap p_j=\varnothing$ 
if $i\ne j$. The set of all such collections will be denoted by
$L(\la/\mu)$. 

With each collection  
$L=(p_i)_{1\le i\le \ell(\la)}\in L(\la/\mu)$ 
we associate a tableau $T$ of 
shape $\la/\mu$ and of order $n$ 
by the following rules: 
$$ 
\gather
\{p_i(t)=(m-1,k),\quad 
p_i(t+1)=(m,k)\}\quad\Rightarrow\quad T(i,m+i)=k;\\ 
\{p_i(t)=(m-1,k-1), \quad
p_i(t+1)=(m,k)\} \quad \Rightarrow\quad T(i,m+i)=k'.
\endgather
$$ 

That is, the $i$th path codes the filling of the $i$th row in
$\la/\mu$, each square $(i,j)$ being associated with a nonvertical
step in $p_i$. If the endpoint of a nonvertical step has coordinates
$(m,k)$ then $j=i+m$, and we have $T(i,j)=k$ or $T(i,j)=k'$ according
to whether the step is horizontal or diagonal. 

To verify that 
$T$ is indeed a tableau of shape $\la/\mu$ 
and of order $n$ we check the four 
properties of Definition A.3. 

Properties \tht{1} and \tht{3}
follow from the definition 
of the path $p_i$. 

Let us check properties \tht2 and \tht4. Given two boxes $(i,j)$
and $(i+1,j)$ in $\la/\mu$, we have to prove that if $T(i,j)$ is
primed then $T(i+1,j)\ge T(i,j)$, and if $T(i,j)$ is nonprimed then
$T(i+1,j)>T(i,j)$. Set $m=j-i$, $k=|T(i,j)|$, $r=|T(i+1,j)|$. The
boxes $(i,j)$ and $(i+1,j)$ correspond to nonvertical steps
$p_i(t)\to p_i(t+1)$ and $p_{i+1}(s)\to p_{i+1}(s+1)$ such that
$p_i(t+1)=(m,k)$ and $p_{i+1}(s+1)=(m-1,r)$, respectively. Since the
paths $p_i$ and $p_{i+1}$ do not intersect, the point $p_{i+1}(s+1)$
must be strictly above the point $p_i(t)$. This condition is exactly
what we need. Indeed, examine the two possible cases: 

First, assume that $T(i,j)$ is primed, $T(i,j)=k'$. This means that
the step $p_i(t)\to p_i(t+1)$ is diagonal, and we have
$p_i(t)=(m-1,k-1)$. Then we get $r>k-1$, i.e., $r\ge k$, as required.

Second, assume that $T(i,j)$ is nonprimed, $T(i,j)=k$. This means that
the step $p_i(t)\to p_i(t+1)$ is horizontal, and we have
$p_i(t)=(m-1,k)$. Then we get $r>k$, as required.

Conversely, if $T$ satisfies the four conditions $(1)-(4)$
then, reversing the above argument, we conclude that the paths are
pairwise nonintersecting. Thus, the 
correspondence $L\rightarrow T$ is a bijection 
between $L(\la/\mu,n)$ 
and $\Tab(\la/\mu,n)$. 

Next, we assign to an arbitrary path $p$ its {\it weight\/} $\Pi(p)$
as follows: 
the weight of $p$ is the product of the weights of its steps, where
 
$\bullet$ a horizontal step with endpoint $(m,k)$ has weight $x_k-a_m$; 
 
$\bullet$ a diagonal step with endpoint $(m,k)$ has weight $y_k+a_m$; 
 
$\bullet$ any vertical step has weight 1.  
  
For an arbitrary collection  
$L=(p_i)_{1\le i\le \ell(\la)}\in L(\la/\mu)$ 
we set 
$$ 
\Pi_{L}=\prod_{i=1}^{\ell(\la)}\Pi(p_i). 
$$ 
Then we get 
$$ 
\Pi_L= 
\prod\Sb \ssq\in\la/\mu \\ T(\ssq)=1,2,\dots,n 
\endSb (x_{T(\ssq)}-a_{c(\ssq)}) 
\prod\Sb \ssq\in\la/\mu \\ T(\ssq)=1',2',\dots,n' 
\endSb (y_{T(\ssq)}+a_{c(\ssq)}), 
$$ 
where $T\leftrightarrow L$.  
It follows that
$$ 
\Sigma_{\la/\mu\mid n}= 
\sum_{L\in L(\la/\mu,n)}\Pi_L.\tag A.4 
$$ 
 
The standard argument
of the Gessel--Viennot theory (see \cite{GV}, \cite{Sa, \S4.5}) 
shows that 
$$ 
\sum_{L\in L(\la/\mu,n)}\Pi_L= 
\det[\Sigma(i,j)]_{1\le i,j\le \ell(\la)},\tag A.5 
$$ 
where $\Sigma(i,j)=\sum_p\Pi(p)$, 
summed over all paths $p$ 
starting at $(\mu_j-j,0)$ and ending 
at $(\la_i-i,n)$. 
{}From (A.3) it follows that 
$$ 
\Sigma(i,j)= 
h_{\la_i-\mu_j+j-i\mid n}(x;y\mid \tau^{\mu_j-j+1}a).\tag A.6 
$$ 
{}From (A.4),(A.5) and (A.6) we obtain that 
$$ 
\Sigma_{\la/\mu\mid n}(x;y\mid a)= 
\det[h_{\la_i-\mu_j+j-i\mid n}
(x;y\mid \tau^{\mu_j-j+1}a)]_{1\le i,j\le \ell(\la)} 
$$ 
and, consequently, 
$$ 
s_{\la/\mu}(x_1,x_2,\dots,x_n;y_1,y_2,\dots,y_n\mid a)= 
\Sigma_{\la/\mu\mid n}(x;y\mid a). 
$$ 
\qed 
\enddemo

\bigskip

{\smc G.~Olshanski}: Dobrushin Mathematics Laboratory, Institute for
Information Transmission Problems, Bolshoy Karetny 19, 
Moscow 101447, GSP-4, Russia.  

E-mail address: {\tt olsh\@iitp.ru, olsh\@online.ru}

{\smc A.~Regev}: Department of Theoretical Mathematics,  
Weizmann Institute of Science, Rehovot 76100, Israel.

E-mail address: {\tt regev\@wisdom.weizmann.ac.il} 

{\smc A.~Vershik}: Steklov Mathematical Institute (POMI), Fontanka 27, 
St.~Petersburg 191011, Russia.

E-mail address: {\tt vershik\@pdmi.ras.ru}

{\smc V.~Ivanov}: Moscow State University.

E-mail address: {\tt vivanov\@vivanov.mccme.ru}

\enddocument 
 
\bye